\documentclass[12pt, a4paper]{amsart}
\usepackage{tikz-cd}
\usepackage{amssymb,amsmath}
\usepackage{array}

\usepackage{amsmath, amssymb, amsfonts, amstext, amsthm, amscd}
\usepackage[mathscr]{euscript}
\usepackage{enumerate}

\usepackage{float}
\usepackage{graphicx}
\usepackage{upgreek}

\usepackage[raiselinks=false,colorlinks=true,citecolor=blue,urlcolor=blue,
linkcolor=blue,bookmarksopen=true,pdftex]{hyperref}

\usepackage{tikz}
\usepackage{tikz-cd}
\usetikzlibrary{snakes}
\usetikzlibrary{intersections, calc}

\usepackage[english]{babel}
\usepackage{chngcntr}

\begin{document}

\thispagestyle{empty}

\def\theequation{\arabic{section}.\arabic{equation}}

\newcommand{\codim}{\mbox{{\rm codim}$\,$}}
\newcommand{\stab}{\mbox{{\rm stab}$\,$}}
\newcommand{\lr}{\mbox{$\longrightarrow$}}

\newcommand{\be}{\begin{equation}}
\newcommand{\ee}{\end{equation}}

\newtheorem{guess}{Theorem}[section]
\newcommand{\bth}{\begin{guess}$\!\!\!${\bf }~}
\newcommand{\eeth}{\end{guess}}
\renewcommand{\bar}{\overline}
\newtheorem{propo}[guess]{Proposition}
\newcommand{\bpropo}{\begin{propo}$\!\!\!${\bf }~}
\newcommand{\epropo}{\end{propo}}

\newtheorem{lema}[guess]{Lemma}
\newcommand{\blem}{\begin{lema}$\!\!\!${\bf }~}
\newcommand{\elem}{\end{lema}}

\newtheorem{defe}[guess]{Definition}
\newcommand{\bdefe}{\begin{defe}$\!\!\!${\bf }~}
\newcommand{\edefe}{\end{defe}}

\newtheorem{coro}[guess]{Corollary}
\newcommand{\bcor}{\begin{coro}$\!\!\!${\bf }~}
\newcommand{\ecor}{\end{coro}}

\newtheorem{rema}[guess]{Remark}
\newcommand{\brem}{\begin{rema}$\!\!\!${\bf }~\rm}
\newcommand{\erem}{\end{rema}}

\newtheorem{exam}[guess]{Example}
\newcommand{\beg}{\begin{exam}$\!\!\!${\bf }~\rm}
\newcommand{\eeg}{\end{exam}}

\newtheorem{notn}[guess]{Notation}
\newcommand{\bnot}{\begin{notn}$\!\!\!${\bf }~\rm}
  \newcommand{\enot}{\end{notn}}

\newcommand{\ch}{{\mathcal H}}
\newcommand{\cf}{{\mathcal F}}
\newcommand{\cd}{{\mathcal D}}
\newcommand{\cR}{{\mathcal R}}
\newcommand{\cv}{{\mathcal V}}
\newcommand{\cn}{{\mathcal N}}
\newcommand{\lra}{\longrightarrow}
\newcommand{\ra}{\rightarrow}
\newcommand{\blr}{\Big \longrightarrow}
\newcommand{\da}{\Big \downarrow}
\newcommand{\ua}{\Big \uparrow}
\newcommand{\hra}{\mbox{{$\hookrightarrow$}}}
\newcommand{\rt}{\mbox{\Large{$\rightarrowtail$}}}
\newcommand{\dua}{\begin{array}[t]{c}
\Big\uparrow \\ [-4mm]
\scriptscriptstyle \wedge \end{array}}
\newcommand{\ctext}[1]{\makebox(0,0){#1}}
\setlength{\unitlength}{0.1mm}
\newcommand{\cl}{{\mathcal L}}
\newcommand{\cp}{{\mathcal P}}
\newcommand{\ci}{{\mathcal I}}
\newcommand{\bz}{\mathbb{Z}}
\newcommand{\cs}{{\mathcal s}}
\newcommand{\ce}{{\mathcal E}}
\newcommand{\ck}{{\mathcal K}}
\newcommand{\cz}{{\mathcal Z}}
\newcommand{\cg}{{\mathcal G}}
\newcommand{\cj}{{\mathcal J}}
\newcommand{\cc}{{\mathcal C}}
\newcommand{\ca}{{\mathcal A}}
\newcommand{\cb}{{\mathcal B}}
\newcommand{\cx}{{\mathcal X}}
\newcommand{\co}{{\mathcal O}}
\newcommand{\bq}{\mathbb{Q}}
\newcommand{\bt}{\mathbb{T}}
\newcommand{\bh}{\mathbb{H}}
\newcommand{\br}{\mathbb{R}}
\newcommand{\bl}{\mathbf{L}}
\newcommand{\wt}{\widetilde}
\newcommand{\im}{{\rm Im}\,}
\newcommand{\bc}{\mathbb{C}}
\newcommand{\bp}{\mathbb{P}}
\newcommand{\ba}{\mathbb{A}}
\newcommand{\spin}{{\rm Spin}\,}
\newcommand{\ds}{\displaystyle}
\newcommand{\tor}{{\rm Tor}\,}
\newcommand{\bff}{{\bf F}}
\newcommand{\bs}{\mathbb{S}}
\def\ns{\mathop{\lr}}
\def\nssup{\mathop{\lr\,sup}}
\def\nsinf{\mathop{\lr\,inf}}
\newcommand{\tG}{{\widetilde{G}}}
\newcommand{\tB}{{\widetilde{B}}}
\newcommand{\tC}{{\widetilde{C}}}
\newcommand{\tW}{{\widetilde{W}}}
\newcommand{\tphi}{{\widetilde{\Phi}}}
\newcommand{\tl}{{\ell}}
\newcommand{\tT}{{\mathtt{t}}}

\title{Equivariant $K$-theory of Springer Varieties}

\author[V. Uma]{Vikraman Uma}
\address{Department of Mathematics, Indian Institute of Technology, Madras, Chennai 600036, India
}
\email{vuma@iitm.ac.in}



\begin{abstract}
  The aim of this paper is to describe the topological equivariant
  $K$-ring, in terms of generators and relations, of a Springer
  variety $\mathcal{F}_{\lambda}$ of type $A$ associated to a
  nilpotent operator having Jordan canonical form whose block sizes
  form a weakly decreasing sequence
  $\lambda=(\lambda_1,\ldots, \lambda_l)$.  This parallels the
  description of the equivariant cohomology ring of
  $\mathcal{F}_{\lambda}$ due to Abe and Horiguchi and generalizes the
  description of ordinary topological $K$-ring of
  $\mathcal{F}_{\lambda}$ due to Sankaran and Uma \cite{su}.
  \end{abstract}

\maketitle

\section{Introduction}\label{intro}
Fix a positive integer $n$ and consider the complete flag
variety $\mathcal F(\mathbb C^n)$ (or more briefly $\mathcal F$)
defined as
\[\mathcal{F}(\mathbb{C}^n):=\{\underline{ V}:=(0=V_0\subset V_1\subset \cdots\subset V_n=\mathbb C^n)\mid
  \dim V_i=i ~~\mbox{for all}~~ i\}.\]

 Let $N:  \mathbb{C}^n\lra \mathbb{C}^n$ 
denote a 
nilpotent linear transformation of $\mathbb C^n$.   
The
Springer variety of type $A$ associated to $N$ denoted by
$\mathcal{F}_{N}$ is the closed subvariety of $\mathcal{F}$
defined as 
\[\{\underline{V}\in \mathcal{F}~~\mid~~NV_i\subset V_{i-1}~~\mbox{for
    all}~~ 1\leq i\leq n\}.\] The Springer variety $\mathcal F_N$ is
seen to be the subvariety of $\mathcal F$ fixed by the action of the
infinite cyclic group generated by the unipotent element
$\mathcal U=I_n+N\in SL(n,\mathbb C)$.  Moreover, denoting by
$\lambda=(\lambda_1,\lambda_2,\ldots,\lambda_l)$ the partition of $n$
where the $\lambda_j$ are the sizes of the diagonal blocks of the
Jordan canonical form of $N$, the variety $\mathcal F_N$ depends, up
to isomorphism, only on the partition $\lambda$. This is so, because
two different choices of nilpotent transformations corresponding to
the same partition $\lambda$ are conjugates in $GL(n,\mathbb C)$.  For
this reason, we assume that $N$ itself is in the Jordan canonical
form:
$N=J_{\lambda}:=\textrm{diag}(J_{\lambda_1},\ldots, J_{\lambda_l})$
with $\lambda_1\ge \cdots\ge \lambda_l$ and denote the Springer
variety $\mathcal F_N$ by $\mathcal F_\lambda$.  (Here
$J_p=(a_{i,j}) \in M_p(\mathbb C)$ is the matrix where
$a_{i,i+1}=1, 1\le i<p,$ and all other entries are zero.)  If
$\lambda=(1,\ldots,1)$, then $N=0$ and we have
$\mathcal F_{\lambda}=\mathcal F(\mathbb{C}^n)=\mathcal F$.  At the
other extreme, when $\lambda=(n)$, $N$ is a regular nilpotent element
and we see that $\mathcal F_{(n)}$ is the one-point variety consisting
only of the standard flag
$0=E_0\subset E_1\subset \cdots\subset E_n=\mathbb C^n$ where $E_j$ is
spanned by the standard basis vectors $e_1,\ldots, e_j$ for
$ 1\le j\le n$.

Note that $\mathcal F_\lambda$ is stable by the action of the
algebraic torus $T^l_\mathbb C\cong (\mathbb C^*)^l$ contained in
$GL(n,\mathbb C)$ consisting of all diagonal matrices which commute
with $N$.  We shall denote by $T^l=(\mathbb S^1)^l$ the compact torus
contained in $T^l_\mathbb C$.  Denoting the diagonal subgroup of
$GL(n,\mathbb C)$ by $T^n_\mathbb C$, we have
$(\tT_1,\ldots, \tT_n)\in T^n_\mathbb C$ belongs to $T^l_\mathbb C$ if
and only if $\tT_{a_j+i}=\tT_{a_j+1}$ for $2\le i\le \lambda_{j+1}$,
$0\leq j\leq l-1$ where
$a_j:=\lambda_1+\cdots+\lambda_j, 1\le j\le l-1$ and $a_0=0$.

The variety $\mathcal F_\lambda$ was first studied by Springer (see
\cite{spr1}, \cite{spr2} and also \cite{hs}).  In particular, Springer
showed that there is a natural action of the symmetric group
$S_n$ on the rational cohomology
$H^*(\mathcal{F}_{\lambda};\mathbb{Q})$ which is compatible with the
standard action of $S_n$ on $H^*(\mathcal{F},\mathbb{Q})$.
Moreover, the restriction homomorphism
$H^*(\mathcal{F};\mathbb{Z})\lra H^*(\mathcal{F}_{\lambda};\mathbb{Z})
$ induced by the inclusion $\mathcal{F}_{\lambda}\hra \mathcal{F}$, is
surjective (see \cite{hs}).  The variety $\mathcal F_\lambda$ is not
irreducible in general, but it is equidimensional.  The irreducible
components of $\mathcal F_\lambda$ are naturally labelled by the set
of standard tableaux of shape $\lambda$.  See \cite{spal}.  Under the
$S_n$-action, the $S_n$-module,
$H^{*} (\mathcal{F}_{\lambda};\mathbb{Q}))$ is
isomorphic to the representation $M_\lambda$ of
$S_n$ induced from the identity representation of the
subgroup
${S}_{\lambda}={S}_{\lambda_1}\times
{S}_{\lambda_2} \cdots \times
{S}_{\lambda_l}\subset {S}_n$.  See \cite{mac} and \cite[p. 204]{dp}.

De Concini and Procesi \cite{dp} gave a description of
$H^*(\mathcal{F}_{\lambda};\mathbb{C})$ as the coordinate ring of an
(unreduced) variety over $\mathbb C$ which we now describe.  Let
$\lambda^\vee$ denote the partition dual to $\lambda$.  The coordinate
ring $\mathbb C[{\mathfrak t}_{\mathbb C}\cap \bar{O}_{\lambda^\vee}]$
of the (non-reduced) scheme
${\mathfrak t}_{\mathbb C}\cap \bar{O}_{\lambda^\vee}$ (scheme
theoretic intersection) where
${\mathfrak t}_{\mathbb C}=Lie(T_\mathbb C^n)\subset
\mathfrak{gl}(n,\mathbb C)=M_n(\mathbb C)$ and
$\bar O_{\lambda^{\vee}}\subset M_n(\mathbb C)$ denotes the closure of
the orbit of $J_{\lambda^\vee}$ under the adjoint action of
$GL(n,\mathbb C)$.  De Concini and Procesi showed that
$H^*(\mathcal F_\lambda;\mathbb C)$ is isomorphic to the algebra
$\mathbb C[{\mathfrak t}_{\mathbb C}\cap \bar O_{\lambda^\vee}]$ (see
\cite{dp}).

Tanisaki \cite{t} described $H^*(\mathcal F_\lambda;\mathbb C)$ as a
quotient of a polynomial ring over $\mathbb{C}$ by an ideal, which has
come to be known as the Tanisaki ideal.  Tanisaki's description in
fact yields the integral cohomology ring of $\mathcal F$.  Recently,
the $T^l$-equivariant cohomology algebra
$H^*_{T^l}(\mathcal F_\lambda;\mathbb Z)$ has been described by H. Abe
and T. Horiguchi.  It turns out that
$H^*_{T^l}(\mathcal F_\lambda;\mathbb Z)$ is the quotient of a
polynomial algebra over $H^*_{T^l}(pt;\mathbb Z)=H^*(BT^l;\mathbb Z)$
modulo an ideal, which is a natural generalization of the Tanisaki
ideal.  This presentation recovers the presentation for the ordinary
integral cohomology ring via the forgetful map
$H^*_{T^l}(\mathcal{F}_{\lambda};\mathbb{Z})\lra
H^*(\mathcal{F}_{\lambda};\mathbb{Z})$ (see \cite[Theorem 4.1]{ah}).

We denote by $\mathcal L_i$ the canonical line bundle over
$\mathcal F(\mathbb C^n)$ whose fibre over a flag
$\underline{V}\in \mathcal{F}({\mathbb C}^n)$ is the vector space
$V_i/V_{i-1}, 1\le i\le n$.  Let
$L_i=\mathcal L_i|_{\mathcal F_\lambda}$.  Recall from \cite{t} that
first Chern classes of $L_i$, $1\leq i\leq n$ generate
$H^*(\mathcal{F}_{\lambda};\mathbb{Z})$.

In \cite{su}, Sankaran and the author described the topological
$K$-ring of $\mathcal{F}_{\lambda}$ in terms of generators and
relations. The generators of $K(\mathcal{F}_{\lambda})$ are the
classes $[L_i]$, $1\leq j\leq n$. The relations are obtained by
interpreting the relations in the cohomology ring in terms of the
classes of the generating line bundles, using suitable gamma
operations in $K$-theory (see \cite[Proposition 4.1]{su}).

It can be seen that $L_i$, $1\leq i\leq n$ are in fact
$T^l$-equivariant line bundles on $\mathcal{F}_{\lambda}$, since they
are the restrictions of the tautological line bundles $\mathcal{L}_i$,
$1\leq i\leq n$ on $\mathcal{F}$ which are $T^n$-equivariant.

In this article we study the $T^l$-equivariant topological $K$-ring of
$\mathcal{F}_{\lambda}$. In our main theorem we give a presentation
for $K^*_{T^l}(\mathcal{F}_{\lambda})$ as an $R(T^l)$-algebra in terms
of generators and relations.

More precisely, we show that $K^0_{T^l}(\mathcal{F}_{\lambda})$ is
generated by the classes $[L_i]_{T^l}$ of the $T^l$-equivariant line
bundles $L_i$ for $1\leq i\leq n$. We further determine the ideal of
relations in $K_{T^l}(\cf_{\lambda})$ as the equivariant analogue of
the $K$-theoretic Tanisaki ideal defined in \cite{su}.

Before stating the main result, we need to set the following
notations.  A non-increasing sequence
$\lambda=(\lambda_1,\ldots, \lambda_l)$ of positive integers where
$\displaystyle\sum_{1\le j\le l} \lambda_j=n$, will be identified with
the partition $(\lambda_1,\ldots, \lambda_n)$ where $\lambda_j=0$ for
$j>l$.  For $1\le s\le n$, let \begin{equation}\label{plambda1}
  p_{\lambda}(s):=\lambda_{n-s+1}+\cdots+\lambda_n.\end{equation}

Recall that
$R(T^l)=K_{T^l}(pt)=\mathbb{Z}[{u_{i}}^{\pm 1}\mid 1\leq i\leq l]$ where
$u_i$, $1\leq i\leq l$ are the characters of $T^l$ corresponding to
the coordinate projections.  

Let $q:=p_{\lambda^{\vee}}(s)$. Let
$\mathcal{R}=R(T^l)[x_1,x_2,\ldots, x_n]$ and let
$\mathcal{I}_{\lambda}$ denote the ideal in $\mathcal{R}$ generated by
the elements \be\label{equivkthrelations}\sum_{0\leq k\leq
  d}(-1)^{d-k}e_{k}(x_{i_1}, x_{i_2},\ldots, x_{i_s})
h_{d-k}({u_{\phi_{\lambda}(1)}},\ldots, {u_{\phi_{\lambda}(s+1-d)}})
\ee for $1\leq s\leq n$, $1\leq i_1<\cdots<i_s\leq n$ and
$d\geq s+1-q$. Here $e_{k}$ stands for the $k$th elementary symmetric
function, $h_{d-k}$ stands for the $(d-k)$th complete symmetric
function (see \cite{mac}) and $\phi_\lambda$ is the map $[n]\lra [l]$
(where $[n]:=\{1,\ldots, n\}$) defined by the condition
\be\label{phi}\begin{array}{ll} (u_{\phi_{\lambda}(1)}, \ldots,
  u_{\phi_{\lambda}(n)}) \\= (\underbrace{u_1,\ldots,
    u_1}_{\lambda_1-\lambda_2},\underbrace{ u_1, u_2,\ldots,
    u_1,u_2}_{2(\lambda_2-\lambda_3)}, \ldots, \underbrace{u_1,\ldots,
    u_l,\ldots, u_1,\ldots,
    u_{l}}_{l(\lambda_l-\lambda_{l+1})})\end{array}\ee as ordered
sequences where for each $1\leq r\leq l$, the $r$th sector of the right
hand side consists of $(u_1,u_2,\ldots, u_r)$ repeated
$(\lambda_r-\lambda_{r+1})$-times. Here we let $\lambda_{l+1}=0$.
(see\cite[(4.2)]{ah}).

We now state our main theorem.

\bth\label{main} With the above notations, let
        \[\Psi_{\lambda}:\mathcal{R}\lra
          K^0_{T^l}(\mathcal{F}_{\lambda})\] be the ring homomorphism
        defined by $\Psi_{\lambda}(x_j)=[L_j]_{T^l}$ for
        $1\leq j\leq n$. Then $\Psi_{\lambda}$ is surjective and
        $\mbox{ker}(\Psi_{\lambda})=\mathcal{I}_{\lambda}$.  \eeth We
        now briefly explain our method of proof of the main theorem
        and an outline of the paper.

Let
\[{n\choose \lambda}:=\frac{n!}{\lambda_1!\cdots \lambda_n!}.\]

Using the fact that $\cf_{\lambda}$ has a $T^l$-stable algebraic cell
decomposition with ${n\choose \lambda}$ cells (see \cite{spal},
\cite{tym}) we show in Theorem \ref{springercell} that
$K^0_{T^l}(\mathcal{F}_{\lambda})$ is a free $R(T^l)$-module of rank
${n\choose \lambda}$. Moreover, we have
$K_{T^l}^1(\mathcal F_\lambda)=0$.

We note that $\mathcal{F}_{\lambda}$ admits a $T^l$-stable filtration
by closed subvarieties arising from the $T^l$-equivariant algebraic
cell decomposition. Using this fact and an induction argument, we show
in Theorem \ref{surj}, that the pull-back map
$K_{T^n}^*(\mathcal{F})\lra K_{T^l}^*(\mathcal{F}_{\lambda})$ is
surjective.

For the special case when $\lambda=(1,1,\ldots, 1)$ the structure of
$K_{T^n}(\cf_\lambda)=K_{T^n}(\cf)$ is well known by the results in
\cite{mcleod} and \cite{kk}. We recall this structure in Section
\ref{equivflag}. In Theorem \ref{eqflagpres} we give a presentation
for $K_{T^n}(\cf)$ as an $R(T^n)$-algebra. In particular, we see that
$K_{T^n}(\cf)$ admits an $S_n$-action.

In Proposition \ref{Snactionspringer} we prove that there exists a natural
$S_n$-action on $K_{T^l}(\cf_{\lambda})$ such that the pull-back map
$K_{T^n}^*(\mathcal{F})\lra K_{T^l}^*(\mathcal{F}_{\lambda})$ is
$S_n$-equivariant. The methods used in proving this result are similar
to those used by Abe and Horiguchi for equivariant cohomology namely,
by restricting to the $T^n$ and $T^l$ fixed points of $\cf$ and
$\cf_{\lambda}$ respectively.

Next, by using the methods similar to those used by Sankaran and the
author for describing the ordinary $K$-ring in \cite[Proposition
4.1]{su} and by the lambda operations in equivariant $K$-theory, we
show that the relations defining the ideal $\mathcal{I}_{\lambda}$
hold in $K_{T^l}^*(\mathcal{F}_{\lambda})$.

This in particular, will show that $\mbox{ker}(\Psi_{\lambda})$
contains $\mathcal{I}_{\lambda}$. Thus we have a well defined
surjective ring homomorphism from $\mathcal{R}/\mathcal{I}_{\lambda}$
to $K_{T^l}^*(\mathcal{F}_{\lambda})$.
        
Now, by using \cite[Lemma 5.2]{ah} we find generators for
$\mathcal{R}/\mathcal{I}_{\lambda}$ as an $R(T^l)$-module.

Thus
$\Psi_{\lambda}:\mathcal{R}/\mathcal{I}_{\lambda}\lra
K_{T^l}(\mathcal{F_{\lambda}})$ is a surjective ring homomorphism
between two free modules of the same rank over $R(T^l)$ which is an
integral domain. Hence it will follow that $\Psi_{\lambda}$ is an
isomorphism. This will prove our main theorem.

Finally, through the forgetful map
$K_{T^l}(\mathcal{F}_{\lambda})\lra K(\mathcal{F}_{\lambda})$ we
recover the presentation of ordinary $K$-ring
$K(\mathcal{F}_{\lambda})$ given in \cite[Theorem 4.2]{su},

\section{Equivariant topological $K$-theory of cellular
  varieties}\label{prelim}

For the definition and basic properties of equivariant topological $K$-ring we
refer to  \cite{segal}.

Let $X$ be a $T_{\mathbb{C}}$-variety for a torus
$T_{\mathbb{C}}\cong (\mathbb{C}^*)^k$ with a $T_{\mathbb{C}}$-stable
algebraic cell decomposition. Let
\[X_{m}\subseteq X_{m-1}\subseteq \cdots\subseteq X_2\subseteq X_1=X\]
be the associated $T_{\mathbb{C}}$-stable filtration so that
$Z_{i}:=X_i\setminus X_{i+1}\cong \mathbb{C}^{k_i}$ for
$1\leq i\leq m$. Here $Z_i\cong \mathbb{C}^{k_i}$ for $1\leq i\leq m$
are the distinct algebraic cells with $T_{\mathbb{C}}$-fixed points
${\mathrm x}_i$ and $\displaystyle X_i=\bigsqcup_{j=i}^m Z_j$. In particular,
$Z_m=X_m=\{{\mathrm x}_m\}$.

We consider the restricted action of the maximal compact subgroup
$T\cong (\mathbb S^1)^k$ of $T_{\mathbb{C}}$ on $X$ as well as on $X_i$ and
  $Z_i$ for $1\leq i\leq m$.

  \bpropo\label{celldec}The ring $K^0_{T}(X)$ is a free $R(T)$-module
  of rank $m$ which is the number of cells. Furthermore, we have
  $K^1_{T}(X)=0$.  \epropo

\begin{proof}
  By \cite[Proposition 2.6, Definition 2.7, Definition 2.8,
  Proposition 3.5]{segal} it follows that we have a long exact
  sequence of $T$-equivariant $K$-groups which is infinite in both
  directions: \be\label{les} \cdots\ra {K}^{-q}_{T}(X_i,X_{i+1})\ra
  {K}_{T}^{-q}(X_i)\ra K^{-q}_{T}(X_{i+1})\ra
  K^{-q+1}_{T}(X_i,X_{i+1})\ra \cdots \ee for $1\leq i\leq m$ and
  $q\in \mathbb{Z}$.

  Moreover, by \cite[Proposition 2.9]{segal} and \cite[Proposition
  3.5]{segal} we have
  \begin{align}K^{-q}_{T}(X_i,X_{i+1}) &=K^{-q}_{T}(X_i\setminus
                                         X_{i+1})\\ &=
                                                      K_{T}^{-q}(\mathbb{C}^{k_i})\cong
                                                      K^{-q}_{T}({\mathrm
                                                      x}_i)\\
                                       &
                                         =\widetilde{K}^{-q}_{T}({\mathrm
                                         x}_i^{+})\\
                                       &=\widetilde{K}_{T}^0({\mathbf
                                         S}^q({\mathrm x}_i^{+}))=R(T)\otimes
                                         \widetilde{K}^0({\mathbf
                                         S}^q({\mathrm x}_i^{+}))\end{align}
                                         for $1\leq i\leq m$(see
                                         \cite[Proposition
                                         2.2]{segal}). Thus when $q$
                                         is even
                                         $K^{-q}_{T}(X_i,X_{i+1})=R(T)$
                                         and when $q$ is odd
                                         $K^{-q}_{T}(X_i,X_{i+1})=0$. Here
                                         ${\mathrm x}_i^{+}$ is the sum of the
                                         $T$-fixed point ${\mathrm x}_i$ and a
                                         base point $\mathfrak{o}$
                                         which is also $T$-fixed (see
                                         \cite[p. 135]{segal}).

                                       Alternately, we can also
                                       identify
                                       $K^{-q}_{T}(X_i,X_{i+1})={\widetilde
                                         K^{-q}_{T}}(X_i/X_{i+1})$
                                       where
                                       $X_i/X_{i+1}\cong
                                       {\mathbf S}^{2k_i}$. For any integer $q$
                                       we have
                                       $\widetilde{K}^{-q}_{T}({\mathbf
                                         S}^{2k_i})\cong
                                       \widetilde{K}_{T}^0({\mathbf S}^{q+2k_i})$
                                       \cite[p.136]{segal}.

                                       Thus when $q$ is even
                                       \[K^{-q}_{T}(X_i,X_{i+1})\cong
                                       \widetilde{K}_{T}^0({\mathbf S}^{q+2k_i})=R(T)\otimes
                                       \widetilde{K}^0({\mathbf S}^{q+2k_i})=R(T)\]
                                       since $q+2k_i$ is
                                       even. Further, when $q$ is odd
                                       \[K^{-q}_{T}(X_i,X_{i+1})=\widetilde{K}_{T}^0({\mathbf
                                           S}^{q+2k_i})=R(T)\otimes
                                       \widetilde{K}^0({\mathbf S}^{q+2k_i})=0\]
                                       since $q+2k_i$ is odd (see
                                       \cite[Section 2, Proposition
                                       3.5, Proposition 2.2]{segal}
                                       and \cite{at}).

                                       Moreover, since
                                       $X_m=Z_m=\{{\mathrm x}_m\}$ and
                                       $X_{m+1}=\emptyset$ we have
                                       $K^0_{T}(X_m)=R(T)$ and
                                       $K_{T}^{-1}(X_m)=K_{T}^{-1}({\mathrm
                                         x}_m)=\widetilde{K}_{T}^{-1}({\mathrm
                                         x}_m^{+})=
                                       K_{T}^0({\mathbf
                                         S}^1({\mathrm x}_m^{+}))=0$ where
                                       ${\mathrm x}_m^{+}={\mathrm x}_m\sqcup
                                       \mathfrak{o}$ where both
                                       ${\mathrm x}_m$
                                       and the base point
                                       $\mathfrak{o}$ are $T$-fixed
                                       (see \cite[p. 135]{segal}).

Now, by decreasing induction on $i$ suppose that
$K_{T}^0(X_{i+1})$ is a free $R(T)$-module for
$1\leq i\leq m$ of rank $m-i$ and $K_{T}^{-1}(X_{i+1})=0$.  We
can start the induction since $K^0_{T}(X_m)=R(T)$ and
$K^{-1}_{T}(X_m)=K_{T}^{-1}({\mathrm x}_m)=0$.

It then follows from \eqref{les} that we get the following split short
exact sequence of $R(T)$-modules \be \label{eq1}
  	\begin{tikzcd}
  		0\arrow{r} & K_{T}^0(X_i, X_{i+1})\arrow{r} & K_{T}^0(X_i)\arrow{r} & K_{T}^0( X_{i+1})\arrow{r} & 0
  	\end{tikzcd}
\ee  
for $1\leq i\leq m$. Thus we have the following
\be\label{eq2} K_{T}^0(X_i)=K_{T}^0(X_{i+1})\bigoplus K_{T}^0(X_i,X_{i+1})\ee

Hence it follows that $K^0_{T}(X_i)$ is a free $R(T)$-module of rank
$m-i+1$. By induction we conclude that $K^0_{T}(X)$ is free
$R(T)$-module of rank $m$ where $X=X_1$.

Since $K^{-1}_{T}(X_i,X_{i+1})=0$ as shown above and
$K^{-1}_{T}(X_{i+1})=0$ by the induction assumption, it also follows
from \eqref{les} that $K^{-1}_{T}(X_i)=0$. Therefore $K_{T}^{-1}(X)=0$
by induction on $i$ since $X=X_1$.
\end{proof}
Henceforth we shall denote $K_{T}^0$ by $K_{T}$.

\section{The equivariant cohomology of
  $\mathcal{F}_{\lambda}$}\label{equivcohom}

\subsection{The $T^n$-equivariant  cohomology of $\mathcal{F}$}
\label{equivcohomfl}

Let $\mathcal V_j$ be the subbundle of the trivial vector bundle
$\mathcal{F}\times \mathbb C^n$ whose fibre over the flag
$\underline{V}=(V_i)\in \mathcal{F}$ is just $V_j$.  Denote by
$\mathcal L_i$ the $T^n$-equivariant line bundle
$\mathcal V_i/\mathcal V_{i-1}, 1\le i\le n,$ on $\mathcal{F}$.

One has an exact sequence of algebraic vector bundles
$0\to \mathcal V_{s-1}\hookrightarrow \mathcal V_s\to \mathcal L_s\to
0$, which leads to an $T^n$-equivariant isomorphism of {\it complex}
vector bundles for $1\le s\le n$:
\begin{equation} \label{vsassumoflinebundles}
\mathcal L_1\oplus \cdots\oplus \mathcal L_s\cong \mathcal V_s.
\end{equation} 
Since, $\mathcal V_n=\epsilon_1\oplus\cdots\oplus \epsilon_n$, the
right hand side is the trivial vector bundle of rank $n$, where the
action of $T^n$ on $\epsilon_i$ is via the character $t_i$
corresponding to the $i$th coordinate projection $T^n\lra S^1$ for
$1\leq i\leq n$. Thus when $s=n$ in \eqref{vsassumoflinebundles}, we
get the following isomorphism of $T^n$-equivariant vector bundles on
$\mathcal{F}$
\begin{equation} \label{sumoftautlinebundles}
\mathcal L_1\oplus \cdots\oplus \mathcal L_n\cong \epsilon_1\oplus\cdots\oplus \epsilon_n.
\end{equation}

Now, by comparing the $T^n$-equivariant $i$th Chern classes in
\eqref{sumoftautlinebundles} we
get \begin{equation}\label{equichern}e_i\left(c_1^{T^n}(\mathcal
    L_1),\ldots, c_1^{T^n}(\mathcal L_n)\right)=e_i(t_1,\ldots,
  t_n)\end{equation} where $e_i$ stands for $i$th elementary symmetric
polynomial and $c_1^{T^n}(\mathcal L_j)$ denotes the $T^n$-equivariant
first Chern class of $\mathcal L_j$ for $1\leq j\leq n$.

Recall that $H^*(BT^n)=\mathbb{Z}[t_1,\ldots, t_n]$. The
map \[H^*(BT^n)[y_1,\ldots, y_n]\lra H_{T^n}^*(\mathcal{F};\mathbb{Z})\]
defined by sending $y_i$ to $c_i^{T^n}(\mathcal L_i)$ induces the following
presentation of $H_{T^n}^*(\mathcal{F};\mathbb{Z})$ as an
$H^*(BT^n)$-algebra. 

\begin{equation}\label{precohomflag} H^*(BT^n)[y_1,\ldots,
  y_n]/\mathcal{J}\cong H^*_{T^n}(\mathcal{F};\mathbb{Z})\end{equation}
where $\mathcal{J}$ is the ideal in $H^*(BT^n)[y_1,\ldots,
y_n]$ generated by the elements
\[e_i(y_1,\ldots, y_n)-e_i(t_1,\ldots, t_n)\] for $1\leq i\leq
n$. (See \cite[Section 3]{ah}.)

\subsection{$T^l$-equivariant cohomology of $\mathcal{F}_{\lambda}$}

The $T^l$-equivariant integral cohomology ring of the Springer variety
$\mathcal F_\lambda$ has been described by Abe and Horiguchi \cite{ah}
in terms of generators and relations, in a way that generalizes the
above description of $H_{T^n}^*(\mathcal F;\mathbb Z)$. We shall
recall below the presentation. We need the following notation (see
\eqref{plambda1}).

\bdefe\label{plambda}
The function $p_\lambda:[n]\to [n]$
associated to a partition $\lambda$ of $n$ defined as
$p_\lambda(s)=\lambda_{n-s+1}+\ldots+\lambda_n,~1\le s\le n$.  
\edefe
Then $p_{\lambda}$ is a monotonically increasing function of $s$.  The
function $p_{\lambda^\vee}$ associated to the {\em dual partition}
$\lambda^\vee$ is more relevant for us.  Recall that the dual
partition $\lambda^\vee $ is defined as
$\lambda^\vee=(\eta_1,\ldots, \eta_n)$ where
$\eta_j=\#\{i\mid \lambda_i\ge j\}$.  Writing $\lambda$ as
$1^{a_1}\cdot 2^{a_2}\cdots n^{a_n}$, where $a_j$ is the number of
times $j$ occurs in $\lambda$, we have $\eta_j=a_j+\cdots+a_n$ for all
$j\ge 1$.

We illustrate in the following example:
\beg Let $n=20$, and,
$\lambda=(5,4,4,2,2,2,1)=1^{1}\cdot 2^{3}\cdot 3^{0}\cdot 4^{2}\cdot 5^{1}$.
Then $\lambda^\vee=(7,6, 3, 3, 1) $ and $p_{\lambda^\vee}(s)=0$ for
$1\le s\le 15$,
$p_{\lambda^\vee}(16)=1,~ p_{\lambda^\vee}(17)=4,~
p_{\lambda^\vee}(18)=7, ~p_{\lambda^\vee}(19)=13,~
p_{\lambda^\vee}(20)=20$.  \eeg

 \bdefe\label{ah ideal} Let $\mathcal{S}=H^*(BT^l)[y_1,\ldots,y_n]$ be
 the polynomial ring in $n$-indeterminates $y_1,\ldots, y_n$ where
 $H^*(BT^l)=\mathbb{Z}[u_1,\ldots, u_l]$ where $u_i$, $1\leq i\leq l$
 are the characters of $T^l$ corresponding to the canonical coordinate
 projections.  The {\em $T^l$-equivariant analogue of the Tanisaki
   ideal} is the ideal $\mathcal{J}_\lambda\subset \mathcal{S}$
 generated by the following elements:
 \[\sum_{k=0}^d (-1)^{d-k}e_k(y_{i_1},\ldots,
   y_{i_s})\cdot h_{d-k}(u_{\phi_{\lambda}(1)},\ldots,
   u_{\phi_{\lambda}(s+1-d)})\] where
 $1\le i_1<\cdots<i_s\le n, ~1\le s\le n$ and 
 $d\ge s+1-p_{\lambda^\vee}(s)$. Here $\phi_{\lambda}$ is as defined
 above \eqref{phi} and $e_k$ and $h_{d-k}$ denote the $k$th
 elementary symmetric function and the $(d-k)$th complete symmetric
 function respectively. (See \cite[(4.1), (4.2)]{ah}) \edefe

 \bth\label{ahtheorem}(\cite[Theorem 4.1]{ah})
Let $\lambda=\lambda_1\ge \cdots\ge \lambda_l$ be a partition of $n$.   Then one has an isomorphism of rings 
\[H_{T^l}^*(\mathcal F_\lambda;\mathbb Z)\cong
  \mathcal{S}/{\mathcal{J}_\lambda }\] where $c^{T^l}_1(L_j)$
corresponds to $y_j+\mathcal{J}_\lambda, 1\le j\le n$.  Moreover, the
inclusion $\iota_\lambda: \mathcal F_\lambda \to \mathcal F$ induces a
surjection
$\iota_\lambda^*: H_{T^n}^*(\mathcal F)\to H_{T^l}^*(\mathcal
F_\lambda)$. Here $c_1^{T^l}(L_j)$ denotes the $T^l$-equivariant first
Chern class of the $T^l$-equivariant line bundle $L_j$. \eeth

The above presentation is the equivariant analogue of the following
presentation of $H^*(\cf_{\lambda},\mathbb{Z})$ due to Tanisaki
\cite{t}.

\brem\label{liftchar} With respect to the inclusion $T^l\hra T^n$
given by sending
\[ (z_1,\ldots, z_l)\mapsto (\underbrace{z_1,\ldots,
    z_1}_{\lambda_1},\underbrace{z_2,\ldots,
    z_2}_{\lambda_2},\ldots,
  \underbrace{z_l,\ldots,z_l}_{\lambda_l}),\] the character $u_i$
is the restriction to $T^l$ of any of the characters $t_j$ of $T^n$
for
$\lambda_1+\cdots+\lambda_{i-1}+1\leq j\leq
\lambda_1+\cdots+\lambda_i$ where $\lambda_1+\cdots+\lambda_{i-1}=0$
when $i=1$.  \erem

\section{$T^n$-equivariant $K$-ring of the flag variety $\mathcal{F}$}\label{equivflag}

When $l=n$ and $\lambda=(1,\ldots, 1)$ then
$\mathcal{F}_{\lambda}=\mathcal{F}$ is the full flag variety.

We recall the notations from Section \ref{equivcohomfl}. We shall
denote the $T^n$-equivariant class of the line bundle $\mathcal L_i$
in $K_{T^n}(\cf)$ by
$[\mathcal{L}_i]_{T^n}\in K_{T^n}^0(\mathcal{F};\mathbb{Z})$ for
$1\le i\le n$. Recall that the $T^n$-equivariant $K$-ring of
$\mathcal{F}$ is well studied classically (see \cite{kk} and
\cite{mcleod}).  Also the presentation for the ordinary $K$-ring of
flag variety and flag bundle are classical and can be found in
\cite[\S3, Chapter IV]{k}.  We give the following presentation for the
$T^n$-equivariant $K$-ring of $\mathcal{F}$ which is analogous to the
presentation of the $K$-ring of the flag bundle as an algebra over the
$K$-ring of the base in \cite{k}.

Let $t_1,\dots, t_n$ denote the characters of $T^n$ corresponding to
the coordinate projections. Thus
$R(T^n)=\mathbb{Z}[t_1^{\pm 1},\ldots, t_n^{\pm 1}]$.
\bth\label{eqflagpres} Let $J'$ be the ideal in
$R(T^n)[x_1,\ldots, x_n]$ generated by the elements
\be\label{flagpres} e_{k}(x_1,\ldots, x_n)-e_k(t_1,\ldots, t_n)\ee for
$1\leq k\leq n$. We have the following isomorphism as a
$R(T^n)$-algebra \be\label{preseqflag} R(T^n)[x_1,\ldots, x_n]/J'
\cong K_{T^n}(\mathcal{F})\ee where $x_i$ maps to
$[\mathcal{L}_i]_{T^n}$ for $1\leq i\leq n$. In particular,
$K_{T^n}(\mathcal{F})$ is generated by $[\mathcal{L}_i]_{T^n}$ for
$1\leq i\leq n$ as an $R(T^n)$-algebra.  \eeth {\bf Proof:} The
unitary group $U(n)\subset GL_n$ acts transitively on the flag variety
$\mathcal{F}$ with stabilizer at the standard flag the diagonal
subgroup $T^n\subset U(n)$. This  identifies $\mathcal{F}$ with the
homogeneous space $U(n)/T^n$.

Further, the line bundle
$\mathcal{L}_i\cong \mathcal{V}_i/\mathcal{V}_{i-1}$ defined above
can be identified with the line bundle
\[U(n)\times _{T^n}\mathbb{C}_i\] on
$\mathcal{F}=U(n)/T^n$ associated to the character
$t_i:T^n\lra S^1$ corresponding to the $i$th coordinate projection for
$1\leq i\leq n$.

Recall from \cite{mcleod} and \cite{kk} that we have an isomorphism
\be\label{mcleodkkiso} R({T}^n)\bigotimes_{R({T}^n)^{S_n}}R({T}^n)\cong
K^0_{{T}^n}(\mathcal F)\ee which is defined on the first factor by
sending $[V]\in R(T^n)$ to the class of the trivial $T^n$-equivariant
vector bundle $[\mathcal{F}\times V]$ on $\mathcal{F}$ (this is
induced by pull back via the constant map $\mathcal{F}\lra pt$) and on
the second factor by sending $[W]\in R(T^n)$ to the class of the
associated vector bundle $U(n)\times_{T^n} W$. In particular, it maps
${t_i}$ to $[\mathcal{L}_i]_{T^n}$ in the second factor. Note that if
we compose the second map with the forgetful homomorphism
$K^0_{T^n}(\mathcal{F})\lra K^0(\mathcal{F})$ we get the classical
Atiyah-Hirzebruch homomorphism $R(T^n)\lra K^0(\mathcal{F})$. Further,
${R({T}^n)^{S_n}}$ is identified with $R(U(n))$ (see \cite{husemoller})
and $R(T^n)$ on the second factor in the tensor product is identified
with $K^0_{U(n)}(\mathcal{F})$.

Recall that $S_n$ acts on
$R(T^n)=\mathbb{Z}[t_i^{\pm1}:1\leq i\leq n]$ by
$\sigma(t_{i}):=t_{\sigma(i)}$ for every $1\leq i\leq n$ and
$\sigma\in S_n$. Thus $R(T^n)^{S_n}$ is
$\mathbb{Z}[t_1,\ldots, t_n]^{S_n}$ localized at the element
$t_1\cdots t_n$, which is invariant under $S_n$. Thus
$R(U(n))=R(T^n)^{S_n}$ is a subring of $R(T^n)$ generated by the elements
$e_{1}(t_1,\ldots, t_n),\ldots, e_n(t_1,\ldots, t_n), 1/(t_1\cdots t_n)$ where
$e_{k}(t_1,\ldots,t_n)$ denotes the $k$th elementary symmetric
polynomial for $0\leq k\leq n$. In particular, note that
$t_1\cdots t_n=e_n(t_1,\ldots, t_n)$ and
$e_0(t_1,\ldots, t_n)=1$ (see \cite{husemoller}).

Since
$\displaystyle R(T^n)\bigotimes_{R(T^n)^{S_n}} R(T^n)\cong
R(T^n)\otimes R(T^n)/J$ where $J$ is the ideal generated by elements
of the form $a\otimes 1-1\otimes a$ where $a\in R(T^n)^{S_n}$. It is
enough the consider $a=e_{k}(t_1,\ldots,t_n)$ and
$a=e_n(t_1,\ldots, t_n)^{-1}$. Thus the ideal $J$ is generated by
\[e_{k}(t_1,\ldots, t_n)\otimes 1-1\otimes e_k(t_1,\ldots, t_n)\] and
\[(t_1^{-1}\cdots t_n^{-1})\otimes 1-1\otimes e_n(t_1,\ldots,
  t_n)^{-1}.\] From the above relations and \eqref{mcleodkkiso} we can
see that in the ring $K^0_{T^n}(\mathcal{F})$ we have the following
\[[\mathcal{L}_i^{\vee}]_{T^n}=(t_1^{-1}\cdots t_n^{-1})\cdot
  \prod_{\substack{j=1 \\j\neq i }}^n[\mathcal{L}_j]_{T^n}\] for
$1\leq i\leq n$. Thus the classes $[\mathcal{L}_i]_{T^n}$ for
$1\leq i\leq n$ generate $K_{T^n}(\mathcal{F})$ as an
$R(T^n)$-algebra.

Consider the $R(T^n)$-algebra homomorphism
$\varphi: R(T^n)[x_1,\ldots, x_n]\ra K_{T^n}(\mathcal{F})$ which
sends $x_i\mapsto [\mathcal{L}_i]_{T^n}$ for $1\leq i\leq n$.


Now, by taking the $k$th exterior power on either side of
\eqref{sumoftautlinebundles} for $1\leq k\leq n$, we get the following
relation in $K_{T^n}(\mathcal{F})$:\be\label{ktheqchern}
e_{k}([\mathcal{L}_1]_{T^n},\ldots,
[\mathcal{L}_n]_{T^n})=e_{k}([\epsilon_1]_{T^n},\ldots,
[\epsilon_n]_{T^n}).\ee

Let $J'$ be the ideal in $R(T^n)[x_1,\ldots, x_n]$ generated by the
elements $e_{k}(x_1,\ldots, x_n)-e_k(t_1,\ldots, t_n)$ for
$1\leq k\leq n$. From \eqref{ktheqchern} it follows that $J'$ is
contained in the kernel of $\varphi$.  Thus $\varphi$ induces a well
defined surjective $R(T^n)$-algebra homomorphism
$\varphi': R(T^n)[x_1,\ldots, x_n]/J'\ra K_{T^n}(\mathcal{F})$.

We know by Proposition \ref{celldec} that $K_{T^n}(\mathcal{F})$ is a
free $R(T^n)$-module of rank $n!$. On the other hand it is well known
that $R(T^n)[x_1,\ldots, x_n]$ is a free
$R(T^n)[e_1(x_1,\ldots, x_n),\ldots, e_n(x_1,\ldots,x_n)]$ module of
rank $n!$ with basis $x_1^{r_1}\cdots x_{n-1}^{r_{n-1}}$ where
$0\leq r_i\leq n-i$ for $1\leq i\leq n-1$ (see \cite[proof of Theorem
3.6, p. 198]{k}). This implies that the quotient
$R(T^n)[x_1,\ldots, x_n]/J'$ is a free $R(T^n)$-module with basis the
products $x_1^{r_1}\cdots x_{n-1}^{r_{n-1}}$. Thus $\varphi'$ is a
well defined and surjective $R(T^n)$-algebra homomorphism between two
free $R(T^n)$-modules of the same rank. Hence $\varphi'$ is an
isomorphism of $R(T^n)$-algebras. Hence the theorem.  $\Box$

\section{$T^l$-equivariant $K$-ring of ${\mathcal F}_{\lambda}$}
In this section we study the $T^l$-equivariant $K$-ring of the
Springer variety ${\mathcal F}_{\lambda}$ with its canonical action of
the torus $T^l$. 

Recall (\cite{segal}) that pull-back via the constant map
$\mathcal{F}_{\lambda}\lra pt$ induces a natural
$R(T^l)=K_{T^l}(pt)$-algebra structure on
$K_{T^l}(\mathcal{F}_{\lambda})$.

\subsection{$R(T^l)$-module structure of
  $K_{T^l}(\mathcal{F}_{\lambda})$}\label{cellular}

\bth\label{springercell} The ring $K^0_{T^l}(\mathcal{F}_{\lambda})$
is a free $R(T^l)$-module of rank $n\choose \lambda$.  In particular,
this implies that $K^0_{T^l}(\mathcal{F}_{\lambda})$ is torsion free
(since $R(T^l)$ is an integral domain). Moreover,
$K_{T^l}^1(\mathcal F_\lambda)=0$.  \eeth \begin{proof} The Springer
  variety $\cf_{\lambda}$ is known to admit a cell decomposition by
  the results of Spaltenstein (see \cite{spal}), with
  ${n\choose \lambda}: =\frac{n!}{({\lambda_1}!\cdots {\lambda_l}!)}$
  locally closed cells isomorphic to affine spaces.  This algebraic
  cell decomposition can be further seen to be $T^l$-stable since the
  cells arise as the intersections of $\cf_{\lambda}$ with the
  $T^n$-invariant Schubert cells in $\cf$(see \cite[Section
  3.2]{tym}). The theorem now follows by Proposition
  \ref{celldec}. \end{proof}

\brem Note that $\mathcal{F}_{\lambda}$ need not be a CW complex but
only has an algebraic cell decomposition (see the example of a string
of pearls in \cite[Section 3.2]{tym}).
\erem

\subsection{The $T^l$-fixed points of $\mathcal{F}_{\lambda}$}\label{fixedpoints}
We recall below  the description of $\mathcal{F}_{\lambda}^{T^l}$ from
\cite[Lemma 2.1]{ah}.

The $T^n$-fixed points of the full flag variety $\mathcal{F}$ are
given by \[\{( \langle e_{w(1)}\rangle \subseteq \langle
  e_{w(1)},e_{w(2)}\rangle\subset \cdots \subset \langle
  e_{w(1)},e_{w(2)},\ldots, e_{w(n)}\rangle=\mathbb{C}^n~\mid~w\in
  S_n\}\]  where $e_1,\ldots, e_n$ is the standard basis of
$\mathbb{C}^n$. Thus we can identify $\mathcal{F}^{T^n}$ with $S_n$.

Then the $T^l$-fixed points of $\mathcal{F}_{\lambda}$ can be
identified with the set of $w\in S_n$ such that $w$ satisfies the following
condition for every $1\leq k\leq l$: \be\label{fp} \mbox{the numbers
  between}~~ \lambda_1+\cdots \lambda_{k-1}+1~~ \mbox{and}~~
\lambda_1+\cdots+\lambda_k\ee appear in the one-line notation of $w$
as a subsequence in the increasing order. We let
$\lambda_1+\cdots+\lambda_{k-1}+1=1$ when $k=1$.

Moreover, $\mathcal{F}_{\lambda}^{T^l}$ can also be identified with
unique right coset representatives of the subgroup $S_{\lambda}$ of
$S_n$ where
$S_{\lambda}:=S_{\lambda_1}\times\cdots \times S_{\lambda_l}$.

Consider the inclusion
$\iota_\lambda:\mathcal F_\lambda \to \mathcal F$ of the Springer
variety in the full flag variety which is $T^l$-invariant with respect
to the restricted $T^l$-action on $\mathcal{F}$ and the natural action
of $T^l$ on $\mathcal{F}_{\lambda}$ described earlier.

This induces a pull back map
\[\iota_\lambda^!:K_{T^n}(\mathcal F)\to K_{T^l}(\mathcal F_\lambda)\]
which factors through $K_{T^l}(\mathcal{F})$ by the inclusion of $T^l$
in $T^n$ ($\iota_{\lambda}^!$ is analogous to the map $\rho_{\lambda}$ in
\cite[Section 3]{ah} for the equivariant cohomology). We can view
$K_{T^l}(\cf_{\lambda})$ as $R(T^n)$-module via the canonical map
$R(T^n)\lra R(T^l)$ induced by the inclusion $T^l\subseteq
T^n$. Moreover the map $R(T^n)\lra R(T^l)$ is surjective since the
character $u_i$ of $T^l$ given by
$(z_1,\ldots, z_{l})\mapsto z_i$ lifts to the character
$t_{\lambda_1+\cdots+\lambda_i}$ of $T^n$ given by
$(\tT_1,\ldots, \tT_n)\mapsto \tT_{\lambda_1+\cdots+\lambda_i}$ for
every $1\leq i\leq l$ (see Remark \ref{liftchar}).

In the following theorem we show that $\iota_\lambda^!$ is
surjective. This statement is analogous in equivariant $K$-theory to
the corresponding statement for ordinary cohomology in \cite{t}, for
equivariant cohomology in \cite{ah} and for ordinary $K$-theory in
\cite{su}.

In the following we let $m={n\choose \lambda}$.


We shall denote by ${\mathrm x}_{w}$ the $T^n$-fixed point in
$\mathcal{F}$ corresponding to $w\in S_n$.



Recall that $\mathcal{F}$ being a nonsingular projective
$T^n$-variety, has a filtrable Bialynicki-Birula cellular
decomposition \cite{bb}. This implies that $\mathcal{F}$ has a
decreasing filtration by $T^n$-stable closed subvarieties:
\be\label{schubertfilt} \{{\mathrm x}_{w_{1}}\}=X_{w_{1}}\subseteq
\cdots \subseteq X_{w_i}\subseteq X_{w_{i+1}}\subseteq \cdots\subseteq
X_{w_{n!}}=\mathcal{F}\ee such that
$X_{w_i}\setminus X_{w_{i-1}}=Z_{w_i}$, where $Z_{w_i}$ is the
Schubert cell corresponding to the $T^n$-fixed point
${\mathrm x}_{w_i}$ for $2\leq i\leq n!$. Moreover,
$X_{w_1}=Z_{w_{1}}=\{{\mathrm x}_{w_{1}}\}$. Here $w_1$ denotes the
identity permutation of $S_n$ and $w_{n!}$ denotes the longest element
with respect to the Bruhat order on $S_n$, where $s_i=(i, i+1)$ for
$1\leq i\leq n-1$ denote the simple reflections in $S_n$.

Let $w_{k_1},\ldots, w_{k_m}$ for $1=k_1<\cdots< k_m\leq n!$ be the
right coset representatives of $S_{\lambda}\backslash S_n$ which
correspond to $\mathcal{F}_{\lambda}^{T^l}$. In particular,
$w_{k_1},\ldots, w_{k_m}$ are the elements of $S_n$ satisfying the
condition \eqref{fp} and $w_{k_1}=w_1$. Furthermore, ${\mathrm x}_{w_{k_j}}$ for
$1\leq j\leq m$ are the $T^l$-fixed points of $\mathcal{F}_{\lambda}$.

\bth\label{surj} The map $\iota_{\lambda}^{!}$ is a surjective
morphism of $R(T^n)$-algebras.  \eeth {\bf Proof:} By \cite{spal}
(also see \cite[Section 3.2]{tym}) there exists a decreasing
filtration of $T^l$-stable closed subvarieties of
$\mathcal{F}_{\lambda}$ given by
\be\label{springerfilt}X'_{w_{k_1}}\subseteq \cdots \subseteq
X'_{w_{k_j}}\subseteq X'_{w_{k_{j+1}}}\subseteq \cdots \subseteq
X'_{w_{k_m}}=\mathcal{F}_{\lambda}\ee such that
${\mathrm x}_{w_{k_j}}\in Z'_{w_{k_j}}=X'_{w_{k_j}}\setminus
X'_{w_{k_{j-1}}}\cong\mathbb{C}^{r_j}$ for $1\leq j\leq m$ where
$X'_{w_{k_{0}}}=\emptyset$. Moreover,
$Z'_{w_{k_j}}=Z_{w_{k_j}}\cap \cf_{\lambda}$. Thus it follows that
$X'_{w_{k_j}}=X_{w_{k_j}}\cap \mathcal{F}_{\lambda}$ for
$1\leq j\leq m$ where $\displaystyle X_{w_{k_j}}=\bigsqcup_{i\leq k_j} Z_{w_i}$ for
$1\leq j\leq m$. Recall that $Z_{w_i}\cap \mathcal{F}_{\lambda}=\emptyset$ for
$i\in [1, n!]\setminus \{k_1,\ldots, k_m\}$.


Consider the following chain of closed subvarieties of $\mathcal{F}$
which is a subchain of \eqref{schubertfilt}:
\[ \{{\mathrm x}_{w_1}\}=X_{w_{k_1}}\subseteq X_{w_{k_{2}}}\subseteq \cdots\subseteq
  X_{w_{k_{m}}}\] where
$\displaystyle X_{w_{k_j}}\setminus X_{w_{k_{j-1}}}=\bigsqcup_{k_{j-1}< i\leq
  k_{j}}Z_{w_i}$ for $1\leq j\leq m$.

Let $\iota_j=\iota_{\lambda}\mid_{X'_{w_{k_j}}}$ denote the inclusion
$X'_{w_{k_j}}\hra X_{w_{k_j}}$ for $1\leq j\leq m$. We shall show that
$\iota_{j}^!:K_{T^n}(X_{w_{k_j}})\lra K_{T^l}(X'_{w_{k_j}})$ is
surjective by induction on $j$.

By following the arguments similar to those in Section \ref{prelim} we
get the following split short exact sequences of $R(T^n)$-modules:
\[0\ra K_{T^n}(X_{w_{k_j}}, X_{w_{k_{j-1}}})\ra
  K_{T^n}(X_{w_{k_j}})\ra K_{T^n}(X_{w_{k_{j-1}}})
  \ra 0 \]
where
$\displaystyle K_{T^n}(X_{w_{k_j}}, X_{w_{k_{j-1}}})\cong
K_{T^n}(X_{w_{k_j}}\setminus X_{w_{k_{j-1}}})=\bigoplus_{k_{j-1}<
  i\leq  k_{j}}K_{T^n}(Z_{w_i})$.
We further have the following split short exact sequence of
$R(T^l)$-modules:
\[0\ra K_{T^l}(X'_{w_{k_j}}, X'_{w_{k_{j-1}}})\ra
  K_{T^l}(X'_{w_{k_j}})\ra K_{T^l}(X'_{w_{k_{j-1}}})
  \ra 0 \] where $K_{T^l}(X'_{w_{k_j}}, X'_{w_{k_{j-1}}})\cong
K_{T^l}(X'_{w_{k_j}}\setminus X'_{w_{k_{j-1}}})= K_{T^l}(Z'_{w_{k_j}})$.

Since $Z'_{w_{k_j}}=Z_{w_{k_j}}\cap \cf_{\lambda}$ is a closed
subvariety of $Z_{w_{k_j}}$ we have the induced map \be\label{res1}
K_{T^n}({\mathrm x}_{w_{k_j}})\cong
K_{T^n}(Z_{w_{k_j}})\stackrel{(\iota_{\lambda}\mid_{Z'_{w_{k_j}
    }})^!}{\lra}
K_{T^l}(Z'_{w_{k_j}})\cong K_{T^l}({\mathrm x}_{w_{k_j}})\ee for every
$1\leq j\leq m$. We note that \eqref{res1} can be identified with the
canonical map $R(T^n)\lra R(T^l)$ and is hence surjective.  It follows
that the map
\be\label{compresproj}\nu_{j}:\bigoplus_{k_{j-1}< i\leq k_{j}}K_{T^n}(Z_{w_i})\ra
K_{T^l}(Z'_{w_{k_j}})\ee obtained by composing the projection to the
factor $K_{T^n}(Z_{w_{k_j}})$ with the map \eqref{res1}, is surjective.

In particular, since
$K_{T^l}(X'_{w_{k_1}})=K_{T^l}(Z'_{w_{k_1}})=K_{T^l}({\mathrm x}_{w_{k_1}})=R(T^l)$,
it follows that the map
$\displaystyle K_{T^n}(X_{w_{k_1}})=K_{T^n}(Z_{w_1})=R(T^n)\ra
R(T^l)= K_{T^l}(X'_{w_{k_1}})$, which can be identified with
the map \eqref{res1} for $j=1$, is surjective.  We assume by induction
that
$\iota_{j-1}^!:K_{T^n}(X_{w_{k_{j-1}}})\ra K_{T^l}(X'_{w_{k_{j-1}}})$
is surjective. Consider the following commuting diagram:

\begin{equation}\label{surjcomm}\begin{array}{llllll} 0 \ra &
                                                              K_{T^n}(X_{w_{k_j}}, X_{w_{k_{j-1}}})\ra &
                                                                                                         K_{T^n}(X_{w_{k_j}})\ra & K_{T^n}(X_{w_{k_{j-1}}})
                                                                                                                                   \ra & 0\\   &~~~~~~~~~~~~\downarrow\nu_j& ~~~~~\downarrow\iota_j^! &~~~~~\downarrow\iota_{j-1}^! & \\ 0\ra & K_{T^l}(X'_{w_{k_j}}, X'_{w_{k_{j-1}}})\ra &
                                                                                                                                                                                                                                                                                        K_{T^l}(X'_{w_{k_j}})\ra & K_{T^l}(X'_{w_{k_{j-1}}})
                                                                                                                                                                                                                                                                                                                   \ra & 0 \end{array}\end{equation}

                                                                                                                                                                                                                                                                                                                 where       the
                                                                                                                                                                                                                                                                                                               rows
                                                                                                                                                                                                                                                                                                               are
                                                                                                                                                                                                                                                                                                               split
                                                                                                                                                                                                                                                                                                               short
                                                                                                                                                                                                                                                                                                               exact
                                                                                                                                                                                                                                                                                                               sequences
                                                                                                                                                                                                                                                                                                               and
                                                                                                                                                                                                                                                                                                               the
                                                                                                                                                                                                                                                                                                               first
                                                                                                                                                                                                                                                                                                               and
                                                                                                                                                                                                                                                                                                               
                                                                                                                                                                                                                                                                                                               third
                                                                                                                                                                                                                                                                                                               vertical
                                                                                                                                                                                                                                                                                                               arrows
                                                                                                                                                                                                                                                                                                               are
                                                                                                                                                                                                                                                                                                               surjections.

                                                                                                                                                                                                                                                                                                               Thus we have 
                                                                                                                                                                                                                                                                                                               $K_{T^n}(X_{w_{k_j}})\cong K_{T^n}(X_{w_{k_j}},X_{w_{k_{j-1}}})\bigoplus K_{T^n}(X_{w_{k_{j-1}}})$ and
                                                                                                                                                                                                                                                                                                               $ K_{T^l}(X'_{w_{k_j}})\cong K_{T^l}(X'_{w_{k_j}},X'_{w_{k_{j-1}}})\bigoplus K_{T^l}(X'_{w_{k_{j-1}}})$. Moreover,
                                                                                                                                                                                                                                                                                                               the
                                                                                                                                                                                                                                                                                                               middle
                                                                                                                                                                                                                                                                                                               vertical
                                                                                                                                                                                                                                                                                                               arrow
                                                                                                                                                                                                                                                                                                               $\iota_j^!$
                                                                                                                                                                                                                                                                                                                  can
                                                                                                                                                                                                                                                                                                                  be
                                                                                                                                                                                                                                                                                                                  identified
                                                                                                                                                                                                                                                                                                                  with the pair $(\nu_j,\iota_{j-1}^!)$ and is therefore 
                                                                                                                                                                                                                                                                                                                  a
                                                                                                                                                                                                                                                                                                                 surjection. This
                                                                                                                                                                                                                                                                                                               completes
                                                                                                                                                                                                                                                                                                               the
                                                                                                                                                                                                                                                                                                               induction. Hence
                                                                                                                                                                                                                                                                                                               $\iota_m^!:
                                                                                                                                                                                                                                                                                                               K_{T^n}(X_{w_{k_m}})\ra
                                                                                                                                                                                                                                                                                                               K_{T^l}(X'_{w_{k_m}})=
                                                                                                                                                                                                                                                                                                               K_{T^l}(\mathcal{F}_{\lambda})$
                                                                                                                                                                                                                                                                                                               is
                                                                                                                                                                                                                                                                                                               surjective.

Now, since $X_{w_{k_m}}$ is a $T^n$-invariant closed subvariety of
$\mathcal{F}=X_{w_{n!}}$, the inclusion $X_{w_{k_m}}\hra \mathcal{F}$ induces
a surjective map of $R(T^n)$-modules
\[K_{T^n}(\mathcal{F})\ra K_{T^n}(X_{w_{k_{m}}}).\]
The theorem now follows by composing the above map with 
$\iota_m^!$.  $\Box$



\subsection{Localization}
When $\lambda=(1,\ldots,1)$ and $\mathcal{F}_{\lambda}=\mathcal{F}$,
it is known (see \cite{kk}) that the map $\iota_1$ in equivariant
$K$-ring induced by restriction to fixed points
$\displaystyle K_{T^n}(\cf)\stackrel{\iota_1}{\lra}
K_{T^n}(\cf^{T^n})=\prod_{w\in S_n} R(T^n) $ is injective.

We have the following result about $\mathcal{F}_{\lambda}$ for an
arbitrary partition $\lambda$.

\blem\label{localization} The canonical restriction
map
\[ K_{T^l}(\cf_{\lambda})\stackrel{\iota_2}{\lra}
  K_{T^l}(\cf_{\lambda}^{T_l})\cong 
  (R(T^l))^m\] is injective where $m={n\choose \lambda}$.
\elem
\begin{proof}
  Since the prime ideal $(0)$ of $R(T^l)$ has support $T^l$ by
  localizing at $(0)$ (see \cite[Proposition 4.1]{segal}), we have
  that
  \[K_{T^l}(\cf_{\lambda})\otimes_{R(T^l)} Q(T^l)\lra
    K_{T^l}(\cf_{\lambda}^{T^l})\otimes_{R(T^l)}Q(T^l)\] is an isomorphism where
  $Q(T^l):=R(T^l)_{\{0\}}$ is the quotient field of the integral
  domain $R(T^l)$. This further implies that the restriction map
  $\displaystyle K_{T^l}(\cf_{\lambda})\lra
  K_{T^l}(\cf_{\lambda}^{T^l})$ is injective since
  $K_{T^l}(\cf_{\lambda})$ is a free $R(T^l)$-module of rank
  $m$ (see Theorem \ref{springercell}).\end{proof}

\subsection{The action of the symmetric group on $K_{T^l}(\mathcal F_\lambda)$}\label{actionS_n}

The following result is the equivariant $K$-theoretic analogue of the
corresponding results for ordinary cohomology (see \cite{spr1},
\cite{spr2} and \cite{hs}), equivariant cohomology (see \cite[Section
3]{ah}) and ordinary $K$-theory (see \cite[Section 3.1]{su}).

\bpropo\label{Snactionspringer} There exists an $S_n$-action on
$K_{T^l}(\cf_{\lambda})$ such
that the map $\iota_{\lambda}^{!}$ is an $S_n$-equivariant
homomorphism.
\epropo

\begin{proof} Our proof is along similar lines as that of the
  corresponding result on equivariant cohomology in \cite[Section
  3]{ah}.  Consider the following commuting square:
  \be\label{commdiag} \begin{array}{llll}
    K_{T^n}(\cf)&\stackrel{\iota_1}{\lra}&
    K_{T^n}(\cf^{T^n})=\bigoplus_{w\in S_n}\mathbb{Z}[t^{\pm
      1}_1,\ldots,
                                           t^{\pm 1}_n]\\
                        \iota_{\lambda}^{!}\da &   & \pi\da \\
                        K_{T^l}(\cf_{\lambda})&\stackrel{\iota_2}\lra
                                         &K_{T^l}(\cf_{\lambda}^{T^l})=\bigoplus_{w\in
                                           \cf_{\lambda}^{T^l}\subseteq
                                           S_n} \mathbb{Z}[u_1^{\pm
                                           1},\ldots, u_l^{\pm 1}]

    \end{array}\ee where the horizontal maps $\iota_1$ and $\iota_2$
  are maps in $T^n$ and $T^l$-equivariant $K$-theory  induced from the
  inclusion of $\cf^{T^n}$ in $\cf$ and $\cf_{\lambda}^{T^l}$ in
  $\cf_{\lambda}$ respectively. As already discussed above the vertical map $\iota_{\lambda}^{!}$
  is induced from the inclusion of $\cf_{\lambda}$ in $\cf$ and is a
  map of $R(T^n)$-algebras. The vertical map $\pi$ is the canonical
  projection $R(T^n)\lra R(T^l)$ induced from the inclusion
  $T^l\subseteq T^n$ on the factors corresponding to $w\in
  \cf_{\lambda}^{T^l}$ and the zero ring map on the other
  factors. 

  We shall now describe the $S_n$-actions on the modules
  $K_{T^n}(\cf)$,\\
  $\displaystyle\bigoplus_{w\in S_n}\mathbb{Z}[t^{\pm 1}_1,\ldots,
  t^{\pm 1}_n]$ and
  $\displaystyle\bigoplus_{w\in \cf_{\lambda}^{T^l}\subseteq S_n}
  \mathbb{Z}[u_1^{\pm 1},\ldots, u_l^{\pm 1}]$. We shall then use
  these to construct an $S_n$ action on $K_{T^l}(\cf_{\lambda})$.

We first recall the left action of the symmetric group $S_n$ on
$K_{T^n}(\cf)$. For this, we consider the right $S_n$-action on the flag
variety $\cf$ as described below (see \cite[Section3]{ah}).

For any $\underline{V}\in \mathcal{F}$ there exists $g\in U(n)$ so
that $\displaystyle V_i=\bigoplus_{j=1}^i \mathbb{C} g(e_j)$, where
$\{e_1,\ldots, e_n\}$ is the standard basis of $\mathbb{C}^n$.  Then
the right action of $w\in S_n$ on $\cf$ can be defined by
\be\label{raflag} \underline{V}\cdot w=\underline{V'}\ee where
$\displaystyle V_i'=\bigoplus_{j=1}^i\mathbb{C} g(e_{w(j)})$. Under
the identification of $\mathcal{F}$ with $U(n)/T^n$ this right action
of $S_n$ corresponds to $gT^n\cdot w:=g\cdot nT^n$ where $n\in N(T^n)$
is a lift of $w$ in $N(T^n)$ which we can choose to be the permutation
matrix associated to $w$ (see \cite[Section 3.11]{kk}).

We recall below the explicit presentation of $K_{T^n}(\cf)$ from Theorem
\ref{eqflagpres}. 

Let $J'$ be the ideal in $R(T^n)[x_1,\ldots, x_n]$ generated by the
elements $e_{k}(x_1,\ldots, x_n)-e_k(t_1,\ldots, t_n)$ for
$1\leq k\leq n$ where
$R(T^n)=\mathbb{Z}[t_1^{\pm 1},\ldots, t_n^{\pm 1}]$. We have the
following isomorphism as an $R(T^n)$-algebra
\[R(T^n)[x_1,\ldots, x_n]/J' \cong K_{T^n}(\mathcal{F})\] where $x_i$
maps to $[\mathcal{L}_i]_{T^n}$ for $1\leq i\leq n$. We shall denote
by $\bar{x}_i$ the class of $x_i$ in $R(T^n)[x_1,\ldots, x_n]/J'$.

The right action \eqref{raflag} of the symmetric group $S_n$ on
$\cf$ induces the following left action on $K_{T^n}(\cf)$:

\be\label{lakflag} w\cdot \bar{x}_{i}:= \bar{x}_{w(i)}~~;~~ w\cdot t_i:=t_i\ee for
$w\in S_n$. This is because the pull back of the line bundle
$\mathcal{L}_i$ under the right action is nothing but the line bundle
$\mathcal{L}_{w(i)}$, and the right action is $T^n$-equivariant.

Now, we shall define a left action of $v\in S_n$ on
$\displaystyle\bigoplus_{w\in S_n}\mathbb{Z}[t_1^{\pm 1},\ldots,
t_n^{\pm 1}]$ as follows: \be\label{symmaction} (v\cdot f)
\mid_{w}=f\mid_{wv}\ee where $f\mid_{w}$ denotes the $w$th component
of $f$, for
$\displaystyle f\in \bigoplus_{w\in S_n}\mathbb{Z}[t_1^{\pm1},\ldots,
t_n^{\pm1}]$ and $w\in S_n$. Note that this is indeed a left action
because
$\left(v_1\cdot (v_2\cdot f)\right)\mid_{w}=(v_2\cdot
f)\mid_{wv_1}=f\mid_{wv_1v_2}=\left((v_1\cdot v_2)\cdot
  f\right)\mid_{w}$.

We note that in the commutative diagram \eqref{commdiag}, the map $\iota_1$ is defined
as follows:
\[\iota_1(\bar{x}_i)\mid_{w}=t_{w(i)}~~;~~\iota_1(t^{\pm 1
  }_i)\mid_{w}=t^{\pm 1}_i\] for $w\in S_n$ and $1\leq i\leq n$.  Thus
it follows that $\iota_1$ is $S_n$ equivariant since
\[\iota_1(w\cdot
  \bar{x}_i)\mid_{v}=\iota_1(\bar{x}_{w(i)})\mid_{v}=t_{v\cdot
    w(i)}=\iota_1(\bar{x}_i)\mid_{v\cdot w}=(w\cdot
  \iota_1(\bar{x}_i))\mid_{v}\] for $1\leq i\leq n$.

Let $\bar{y}_i$ denote the image in $K_{T^l}(\cf_{\lambda})$ of
$\bar{x}_i$ under the surjective $R(T^n)$-algebra homomorphism
$\iota_{\lambda}^{!}$. We have the following lemma which is analogous
to \cite[Lemma 3.1]{ah}. This follows from Theorem \ref{surj} and
Theorem \ref{eqflagpres}.

\blem\label{geneqkspringer} The $T^l$-equivariant topological $K$-ring
$K_{T^l}(\cf_{\lambda})$ is generated by $\bar{y}_1,\ldots, \bar{y}_n$
as an algebra over
$R(T^l)=\mathbb{Z}[u_1^{\pm 1},\ldots, u_l^{\pm 1}]$.  \elem Now we
proceed to construct an $S_n$-action on
$\displaystyle\bigoplus_{w\in \cf_{\lambda}^{T^l}}\mathbb{Z}[u^{\pm
  1}_1,\ldots, u^{\pm 1}_l]$ and on $K_{T^l}(\cf_{\lambda})$.  Let
$\phi:[n]\longrightarrow [l]$ be the map defined
by \begin{equation}\label{phi_1}\phi(i)=k\end{equation} if
$\lambda_1+\cdots+\lambda_{k-1}+1\leq i\leq
\lambda_1+\cdots+\lambda_k$ where $\lambda_1+\cdots+\lambda_{k-1}=0$
when $k=1$ (see \cite[p.6]{ah}). We note that the inclusion
$T^l\hra T^n$ induces the map $R(T^n)\lra R(T^l)$ given by
$t_i^{\pm1}\mapsto u_{\phi(i)}^{\pm1}$ for $1\leq i\leq n$.
  
  Then the map  
  \[\pi:\bigoplus_{w\in S_n} \mathbb{Z}[t_1^{\pm
      1},\ldots, t_n^{\pm 1}]\lra \bigoplus_{w\in
      \cf_{\lambda}^{T^l}}\mathbb{Z}[u^{\pm 1}_1,\ldots, u^{\pm
      1}_l]\] in the commutative diagram \eqref{commdiag} is given by
\[ \pi(f\mid_{w}(t^{\pm 1}_1,\ldots, t^{\pm 1}_n))=f\mid_{w}(u^{\pm
    1}_{\phi(1)},\ldots, u^{\pm 1}_{\phi(n)}).\]

Thus from the commutative diagram \eqref{commdiag}, for
$w\in \mathcal{F}_{\lambda}^{T^l}$ we have
\be\label{r1}\iota_2(\bar{y}_i)\mid_w=\iota_2(\iota_{\lambda}^{!}(\bar{x_i}))\mid_{w}=\pi(\iota_1(\bar{x}_i)\mid_{w})=\pi(t_{w(i)})=u_{\phi(w(i))}\ee
and \be\label{r2}\iota_2(u_i^{\pm 1})\mid_w=u_i^{\pm 1}.\ee

The left action of $v\in S_n$ on
$\displaystyle\bigoplus_{w\in \cf_{\lambda}^{T^l}}\mathbb{Z}[u_1^{\pm
  1},\ldots, u_l^{\pm 1}]$ is defined by \be\label{r3}(v\cdot
f)\mid_w=f_{w'}\ee for $w\in \cf_{\lambda}^{T^l}$ and
$\displaystyle f\in \bigoplus_{w\in
  \cf_{\lambda}^{T^l}}\mathbb{Z}[u_1^{\pm 1},\ldots, u_l^{\pm 1}]$ and
$w'\in \cf_{\lambda}^{T^l}$ is the coset representative of the right
coset $[wv]$. (Recall that $\cf_{\lambda}^{T^l}$ can be identified
with the coset representatives of the cosets of the subgroup
$S_{\lambda}$ in $S_n$.)

We have the following lemma analogous to \cite[Lemma 3.2]{ah}.
\blem\label{equiviota2} For every $v\in S_n$, $1\leq i\leq n$ and
$1\leq j\leq l$ we have \be\label{r4} v\cdot
(\iota_2(\bar{y}_i))=\iota_2(\bar{y}_{v(i)}) \ee and \be\label{r5}
v\cdot (\iota_2(u_j^{\pm 1}))=\iota_2(u_j^{\pm 1}).\ee \elem
\begin{proof}
  From \eqref{r2} and \eqref{r3} we have
  \[(v\cdot
    \iota_2(u^{\pm1}_i))\mid_{w}=\iota_{2}(u^{\pm1}_i)\mid_{w'}=u^{\pm1}_i=\iota_{2}(u^{\pm1}_i)\mid_{w}\]
  for every $w\in S_n$. Thus \eqref{r5} follows.

  Now, from \eqref{r1} and \eqref{r3} we have
  \[v\cdot
    \iota_2(\bar{y}_i)\mid_{w}=\iota_2(\bar{y}_i)\mid_{w'}=u_{\phi(w'(i))}\]
  and
  \[\iota_2(\bar{y}_{v(i)})\mid_{w}=u_{\phi(w(v(i)))}.\]
  Thus it suffices to prove that $\phi(w'(i))=\phi(wv(i))$. This
  follows since $[w']=[wv]$ in $S_{\lambda}\backslash S_n$ (see
  \cite[proof of Lemma 3.2]{ah}). Hence the proof.
\end{proof}  
Since $\iota_2$ is injective we therefore obtain an $S_n$ action on
$K_{T^l}(\cf_{\lambda})$ given by \be\label{actionspringerkring}
w\cdot \bar{y}_i=\bar{y}_{w(i)} ~\mbox{and} ~~w\cdot u_j^{\pm
  1}=u_j^{\pm 1}\ee for $w\in
S_n$ for $1\leq i\leq n$ and $1\leq j\leq l$. It follows that the action is well defined by Lemma
\ref{geneqkspringer} and Lemma \ref{equiviota2}.

By identifying $\bar{y}_i$ with $[L_i]_{T^l}$ the $S_n$-action on
$K_{T^l}(\cf_{\lambda})$ is given by
$w\cdot [L_i]_{T^l}:=[L_{w(i)}]_{T^l}$ for $1\leq i\leq n$ and
$w\cdot u^{\pm1}_j=u^{\pm1}_j$ for $1\leq j\leq l$.

Moreover, this also
implies that for the $S_n$ action on $K_{T^n}(\cf)$ and
$K_{T^l}(\cf_{\lambda})$ the map $\iota_{\lambda}^{!}$ is
$S_n$-equivariant.
\end{proof}

\subsection{Sectioning canonical bundles over $\mathcal F_\lambda$}

Recall the function $p_\lambda$ defined in Definition \ref{plambda}.
The numbers $p_{\lambda^\vee}(s)$ are related to the nilpotent
transformation $N=J_{\lambda}$ as follows.

\blem\label{rank} \cite[Proposition 3]{t}.
With the above notations, \\(i) $p_{\lambda^\vee}(s)=\mathrm{rank}(J_{\lambda}^{n-s}), ~1\le s\le n$.\\
(ii) Let
$\underline V=V_1\subset V_2\subset \cdots\subset V_n=\mathbb C^n$ be
a flag that refines the partial flag
$0=\im (J_{\lambda}^{\lambda_1})\subset
\im(J_{\lambda}^{\lambda_1-1})\subset\cdots \subset \im
(J_{\lambda}^2)\subset \im(J_{\lambda})\subset \mathbb C^n$.  Then,
for any $\underline V'\in \mathcal F_\lambda$ and any $s\ge 1$, we
have $V_q\subset V'_s$ where $q=p_{\lambda^\vee}(s)$.  \elem

\bpropo \label{relationsinflambda} For $1\le s\le n$, we have the
isomorphism \be\label{isom} L_{1}\oplus \cdots \oplus L_{s}\cong
\xi\oplus \epsilon_{j_1}\oplus\cdots \oplus\epsilon_{j_q}\ee for some
$T^l$-equivariant complex vector bundle $\xi=\xi(\mathbf{i})$ of rank
$s-q$ and $T^l$-equivariant trivial line bundles $\epsilon_{j_{\tl}}$
for $1\leq {\tl}\leq q$ over $\mathcal F_\lambda$ where
$q:=p_{\lambda^\vee}(s)$.  Moreover, for $d\geq s+1-q$, \eqref{isom}
also implies the following isomorphism \be\label{isom1} L_{1}\oplus
\cdots \oplus L_{s}\cong \xi'\oplus \epsilon_{j_1}\oplus\cdots
\oplus\epsilon_{j_{s+1-d}}\ee where
\[\xi'=\xi\oplus \epsilon_{j_{s+1-d+1}}\oplus \cdots \oplus
  \epsilon_{j_q}\] is a $T^l$-equivariant complex vector bundle of
rank $d-1$. Furthermore, we can choose $\xi$ and $\epsilon_{j_{\tl}}$
$1\leq \tl\leq q$ such that the $T^l$ action on the fibre
$\mathbb{C}^q$ of the trivial bundle
$\epsilon_{j_1}\oplus\cdots \oplus\epsilon_{j_q}$ of rank $q$ is via
the tuple of characters
$(u_{\phi_{\lambda}(1)},\ldots, u_{\phi_{\lambda}(q)})$ \epropo

\begin{proof}
  We replace $N$ by a conjugate $gN g^{-1}$ so that
  $\im (gNg^{-1})^{n-k}=E_{p_{\lambda^\vee}(k)}=\mathbb
  C^{p_{\lambda^\vee}(k)}$ for $k\ge 1$.  We may then choose the
  refinement of
  $(\cdots\subset (gNg^{-1})^2\mathbb{C}^n\subset
  (gNg^{-1})\mathbb{C}^n\subset \mathbb{C}^n)$ to be the standard flag
  $\underline{E}=(0=E_0\subset E_1\subset \cdots\subset E_n=\mathbb{C}^n)\in
  \mathcal F$ where $E_i=\langle e_1,\ldots, e_i\rangle$ for
  $1\leq i\leq n$.  Thus by Lemma \ref{rank}, we have
  $\mathbb C^q\subset V_s$ for any
  $\underline{V}\in \mathcal F_{gNg^{-1}}$.  Let
  $\iota_g:\mathcal F\to \mathcal F$ be the translation by $g$:
  $\underline {V}\mapsto g\underline {V}=gV_0=0\subset gV_1\subset
  \cdots \subset gV_n=\mathbb C^n$.  Since $GL(n,\mathbb C)$ is
  connected, the composition
  $\mathcal F_{N}\stackrel{\iota_\lambda}{\hookrightarrow} \mathcal
  F\stackrel{\iota_g}{\longrightarrow }\mathcal F$, denoted
  $\iota_{\lambda,g}$ is homotopic to $\iota_\lambda$ and maps
  $\mathcal F_{N}$ onto $\mathcal F_{gNg^{-1}}\subset \mathcal F$.  It
  follows that
  $\iota^{!}_{\lambda,g}: K_{T^n}(\mathcal{F})\lra
  K_{T^l}(\mathcal{F}_{\lambda})$ can be identified with
  $\iota^{!}_\lambda$.  In particular,
  $\iota_{\lambda,g}^!([\mathcal L_j])=[L_j]~\forall ~1\leq j\le n$.

  Recall that the left action of $g\in GL_n(\mathbb{C})$ on
  $(x_1,\ldots, x_n)\in\mathbb{C}^n$ is defined by multiplying the
  $n\times n$ matrix corresponding to $g$ from the left on the column
  vector $(x_1,x_2,\ldots, x_n)^{\mathrm t}$.  Note that for a
  representative $g$ of $w\in S_n$ in $GL_n(\mathbb{C})$, the left
  action $\iota_g$ of $g$ on $\mathcal{F}$, which is induced by the
  left action of $w$ on $\mathbb{C}^n$ given by
  $w\cdot (x_1,\ldots,x_n)=(x_{w^{-1}(1)},\ldots, x_{w^{-1}(n)})$,
  satisfies $\iota_{g}^*(\mathcal L_i)\simeq \mathcal L_i$, for
  $1\leq i\leq n$. Moreover, $\iota_{g}$ is $T^n$ equivariant with
  respect to a group homomorphism $\psi_g:T^n\lra T^n$ given by
  \be\label{equivaction}(\tT_1,\ldots, \tT_n)\mapsto
  (\tT_{w^{-1}(1)},\ldots, \tT_{w^{-1}(n)}).\ee Note that $\psi_g$ is
  conjugation by $g$. In particular, under $\psi_g$, an element
  $(z_1,\ldots, z_l)\in T^l$ maps to
  $(z_{\phi(w^{-1}(1))},\ldots, z_{\phi(w^{-1}(n))})$ (see Remark
  \ref{liftchar} and \eqref{phi_1}). Thus the image of $T^l$ under $\psi_g$ is the
  conjugate subgroup $gT^lg^{-1}$ of $T^n$.  We shall consider the
  action of $T^l$ on $\mathcal{F}$ and $\mathcal{F}_{gNg^{-1}}$ via
  $\psi_g$.

  We remark here that under the isomorphism of $\mathcal{F}$ with
  $U(n)/T^n$ this left action of $w\in S_n$ corresponds to
  $w\cdot gT^n=ngT^n$, where $n\in N(T^n)$ can be chosen to be the
  permutation matrix associated to $w$.

Recall from \cite[Section 5]{ah} that there exists
$\bar{w}\in \mathcal{F}_{\lambda}^{T^l}\subseteq S_n$ such that
$\bar{w}\cdot \underline{E}$ refines the flag
$(\cdots\subset N^2\mathbb{C}^n\subset N\mathbb{C}^n\subset
\mathbb{C}^n)$. Here $\underline{E}=(E_i)\in \mathcal{F}$
denotes the standard flag with $E_i=\langle e_1,\ldots, e_i\rangle$
for $1\leq i\leq n$. In particular, we have
$\bar{w}\cdot E_q=\bar{w}\cdot \mathbb{C}^q=\mbox{Im}(N^{n-s})$ where
$q=p_{\lambda^{\vee}}(s)$.

Thus conjugating $N$ by a representative
$g\in GL_n(\mathbb{C})$ of $\bar{w}^{-1}$ we get
\[\mbox{Im}(gNg^{-1})^{n-s}=E_q=\mathbb{C}^q\] where
$q=p_{\lambda^{\vee}}(s)$. This implies that $\underline{E}$ will
refine the flag
\[(\cdots\subset (gNg^{-1})^2\mathbb{C}^n\subset
(gNg^{-1})\mathbb{C}^n\subset \mathbb{C}^n).\]

We further recall from \cite{ah} that
\be\label{relphilambda}\phi_{\lambda}=\phi\circ \bar{w}\ee where $\phi_{\lambda}$ is as
defined in \eqref{phi} and $\phi$ is as in \eqref{phi_1}).  

With this choice of $\bar{w}$ and $g$ and the preceeding arguments, it
follows that the image of $(z_1,\ldots, z_{l})\in T^l$
under $\psi_g$ is
$(z_{\phi\circ \bar{w}(1)},\ldots, z_{\phi\circ \bar{w}(n)})$. Thus by
\eqref{relphilambda} it follows that the $T^l$-action on
$\mathbb{C}^n$ is via the tuple of characters
$(u_{\phi_{\lambda}(1)},\ldots, u_{\phi_{\lambda}(n)})$.

With the above assumptions, for the rest of the proof we shall denote
$gNg^{-1}$ simply by $N$ and the image $\psi_g(T^l)=gT^lg^{-1}$ simply
by $T^l$.

Let $G_{n,s}=G_s(\mathbb C^n)$ denote the Grassmann variety of
$s$-planes in $\mathbb C^n$.  One has a projection
$\pi_s: \mathcal F\to G_{n,s}$ defined as $\underline{V}\mapsto
V_s$. Note that the map $\pi_s$ is $T^l$-equivariant where the $T^l$
action on $G_{n,s}$ is the natural one obtained by restricting the
action of $T^n$ which takes an $s$-dimensional subspace of
$\mathbb{C}^n$ to another $s$-dimensional subspace.

Let $Y_q\subset G_{n,s} $ denote the subvariety
$\{U\in G_{n,s}\mid U\supset \mathbb{C}^q\}, 1\le q<s$. Note that
since $T^l$ commutes with $N$ it stabilizes
$Im(N)^{n-s}=E_q=\mathbb{C}^q$. Thus the inclusion $Y_q\hra G_{n,s} $
is $T^l$-equivariant.

Moreover, $Y_q$ is isomorphic to a Grassmann variety $G_{n-q,s-q}$.  A
specific isomorphism $Y_q\cong G_{s-q}(\mathbb C^n/\mathbb{C}^q)$ is
obtained by sending $U\in Y_q$ to $U/\mathbb{C}^q$.  The tautological
complex vector bundle $\gamma_{n,s}$ is of rank $s$, whose fibre over
$A\in G_{n,s}$ is the vector space $A$. Note that $\gamma_{n,s}$ is
$T^l$-equivariant vector bundle. When restrictred to $Y_q$,
$\gamma_{n,s}$ has a trivial subbundle $q\epsilon $ of rank $q$.  The
subbundle $q\epsilon$ of $\gamma_{n,s}$ is also $T^l$-equivariant
since the action of $T^l$ preserves the subspace $\mathbb{C}^q$.  We note that
$T^l$ acts on $\mathbb{C}^q\subseteq \mathbb{C}^n$ via the tuple of
characters $(u_{\phi_{\lambda}(1)},\ldots, u_{\phi_{\lambda}(q)})$.

Indeed we have a commuting diagram
\[\begin{array}{ccc}Y_q\times \mathbb{C}^q& \hookrightarrow&  E(\gamma_{n,s}|_{Y_q})\\
    \downarrow && \downarrow \\
    Y_q&\stackrel{id}{\longrightarrow}& Y_q\\
  \end{array}\]
where all maps are $T^l$-equivariant and the vertical arrows are
bundle projections. Thus as a $T^l$-equivariant vector bundle we shall write
$Y_q\times \mathbb{C}^q$ as
$\epsilon_{\chi_1}\oplus\cdots \oplus \epsilon_{\chi_q}$.
Therefore
\begin{equation} \label{splitting} \gamma_{n,s}|_{Y_q}\cong
  \omega\oplus \epsilon_{\chi_1}\oplus\cdots \oplus \epsilon_{\chi_q}
\end{equation}
where $\omega$ is the complex vector bundle over $Y_q$ whose fibre
over $A\in Y_q$ is the complex vector space
$A':=A/\mathbb{C}^q$. Since $\gamma_{n,s}$ and
$\epsilon_{\chi_1}\oplus\cdots \oplus \epsilon_{\chi_q}$ are
$T^l$-equivariant we get that the quotient bundle $\omega$ is
$T^l$-equivariant. Thus the direct sum decomposition \eqref{splitting}
is $T^l$-equivariant.


It follows by Lemma \ref{rank} that the image of the
composition
\[\mathcal{F}_{N}\stackrel{\iota_{N}}{\hookrightarrow}\mathcal F\stackrel{\pi_s}{\longrightarrow} G_{n,s},\] denoted
$\pi_{N,s}$, is contained in $Y_q$.  Furthermore, the map
$\pi_{N,s}$ is $T^l$-equivariant.  Thus we have a
commuting diagram where all the maps are $T^l$-equivariant
\begin{equation}
\begin{array}{rcc}
\mathcal F_{N} &\stackrel{\iota_{N}}{\hookrightarrow} & \mathcal F\\
\pi_{N,s}\downarrow &&~~\downarrow \pi_s\\
Y_q& \hookrightarrow & G_{n,s}.\\
\end{array}
\end{equation}
Now,
\[\displaystyle\pi_s^*(\gamma_{n,s})=\mathcal V_s=\mathcal L_1\oplus \cdots \oplus
  \mathcal L_s\] by (\ref{vsassumoflinebundles}).  Also recall that
$\iota_{N}^*(\mathcal{L}_i)=L_i$ for $1\leq i\leq n$.

Therefore
\begin{align} L_1\oplus \cdots\oplus L_s&=\iota_{N}^*(\mathcal
  L_1\oplus\cdots\oplus \mathcal L_s)\\ &=\iota_{N}^*\circ
  \pi_s^*(\gamma_{n,s})\\ &=\pi_{N,s}^*(\gamma_{n,s}|_{Y_q})\\
  &= \pi_{N,s}^*(\omega)\oplus
  \pi_{N,s}^*(\epsilon_{\chi_1}\oplus\cdots \oplus
  \epsilon_{\chi_q}),\end{align} from (\ref{splitting}). Thus if
$\epsilon_{j_{\tl}}:=\pi_{N,s}^*(\epsilon_{\chi_{\tl}})$ for
$1\leq \tl\leq q$ and $\xi:=\pi_{N,s}^* (\omega)$ then the
isomorphism \eqref{isom} holds over $\mathcal{F}_{N}$.

Moreover, in the isomorphism \eqref{isom} over $\mathcal{F}_{N}$, $T^l$
acts on the fibre of the $T^l$-equivariant trivial bundle of rank $q$
\[\bigoplus_{\tl=1}^q\epsilon_{j_{\tl}}\cong
  \mathcal{F}_{N}\times \mathbb{C}^q\] via the tuple
$(u_{\phi_{\lambda}({1})}, \ldots, u_{\phi_{\lambda}(q)})$. Hence the
proposition.




\end{proof}

\subsection{The lambda operations in equivariant $K$-theory.}  We
recall here the lambda operations or the exterior power operations
in equivariant $K$-theory (see \cite[Section 1]{atiyahtall}).

Let $X$ be a compact $G$-space for $G$ a compact Lie group. For
$x\in K_{G}(X)$ we let
\[\uplambda_t(x):=\sum_{i\ge 0} \uplambda^i (x)t^i\] as an element of the
formal power series ring $K_{G}(X)[[t]]$ in the indeterminate $t$.
For any $x\in K_{G}(X)$, $\uplambda^0(x)=1$. When
$x=[\xi]_{G}\in K_{G}(X)$ is the class of a $G$-vector bundle $\xi$ on
$X$ of rank $k$, $\uplambda_t(x)$ is a polynomial of degree $k$, since
the exterior power $\uplambda^d([\xi]_{G})=[\Lambda^d(\xi)]_{G}=0$ for
$d\ge k+1$.  In particular, when $\xi$ is a $G$-line bundle, then
$\uplambda_t([\xi]_{G})=1+[\xi]_{G}t$.

Let $\xi=\xi_{1}\oplus\cdots\oplus \xi_{k}$ where $\xi_i$ is a
$G$-equivariant line bundle on $X$ for $1\leq i\leq k$. Then we have
$\displaystyle\uplambda_t([\xi]_{G})=\prod_{i=1}^k(1+[\xi_i]_{G}t)$ and
$\displaystyle\uplambda_t([-\xi]_{G})=\prod_{i=1}^k(1+[\xi_i]_{G}t)^{-1}$. The
last two equalities follow from the identity
\[\uplambda_t(x+y)=\uplambda_t(x)\cdot\uplambda_t(y).\]

When $X=\cf_{\lambda}$ and $G=T^l$ then from \eqref{isom1} and the
above properties of $\uplambda_t$ we have
\be\label{lambdadirectsum1}\uplambda_t([\xi']_{T^l})=\prod_{1\leq r\leq
  s}(1+[L_{r}]_{T^l}t)\prod_{1\leq \tl\leq
  {s+1-d}}(1+[\epsilon_{j_{\tl}}]_{T^l}t)^{-1}\ee for $d\geq s+1-q$
where $q=p_{\lambda^{\vee}}(s)$. Now, $\uplambda_t([\xi']_{T^l})$ is a
polynomial of degree $d-1$ since $\mbox{rank}(\xi')=d-1$ for
$d\geq s+1-q$.  Thus by comparing the coefficients of $t^d$ in
\eqref{lambdadirectsum1} where $d\geq s-q+1$ we obtain the following
equation in $K_{T^l}(\cf_{\lambda})$ \be\label{relations} \sum_{0\leq
  k\leq d} (-1)^{d-k}e_{k}([L_{1}]_{T^l},\ldots,
[L_{s}]_{T^l})h_{d-k}([\epsilon_{j_1}]_{T^l},\ldots,
[\epsilon_{j_{s+1-d}}]_{T^l})=0.\ee We have the following proposition
with the notations of Section \ref{intro} where
$\mathcal{I}_{\lambda}$ in particular denotes the equivariant
$K$-theoretic Tanisaki ideal (see \eqref{equivkthrelations}).
\bpropo\label{wdsurj} We have a well defined surjective
$R(T^l)$-algebra homomorphism
  \[\bar{\Psi_{\lambda}}:R(T^l)[x_1,\ldots, x_n]/\mathcal{I}_{\lambda}\lra
    K_{T^l}(\cf_{\lambda})\] which sends $x_i$ to
  $[L_i]_{T^l}\in K_{T^l}(\cf_{\lambda})$ for $1\leq i\leq n$. \epropo
  \begin{proof}
    Recall from Theorem \ref{surj} that $K^{T^l}(\cf_{\lambda})$ is
    generated by $[L_i]_{T^l}$ as an $R(T^l)$-algebra. Furthermore,
    recall from Proposition \ref{Snactionspringer} that the $S_n$ action on
    $K_{T^l}(\cf_{\lambda})$ is given by
    $w\cdot [L_{i}]_{T^l}:=[L_{w(i)}]_{T^l}$ and
    $w\cdot u_i^{\pm 1}=u_i^{\pm 1}$ where
    $R(T^l)=\mathbb{Z}[u_1^{\pm1},\ldots, u_l^{\pm1}]$.  By
    Proposition \ref{relationsinflambda}, for $1\leq \tl\leq q$,
    $\epsilon_{j_{\tl}}$ is a trivial line bundle on
    $\mathcal{F}_{\lambda}$ with $T^l$ action on the fibre given by
    the character $u_{\phi_{\lambda}(\tl)}$. Thus we can identify
    $[\epsilon_{j_{\tl}}]_{T^l}$ with the element
    $u_{\phi_{\lambda}(\tl)}\in R(T^l)$ for $1\leq \tl\leq q$.  Now,
    through the $S_n$-action the relation \eqref{relations} implies
    the set of relations \be\label{relations1} \sum_{0\leq k\leq d}
    (-1)^{d-k}e_{k}([L_{i_1}]_{T^l},\ldots,
    [L_{i_s}]_{T^l})h_{d-k}(u_{\phi_\lambda(1)},\ldots,
    u_{\phi_{\lambda}(q)})=0\ee in $K_{T^l}(\mathcal{F}_{\lambda})$
    for $1\leq i_1<\cdots<i_s\leq n$, $1\leq s\leq n$ and
    $d\geq s-q+1$. Hence the proposition.
    \end{proof}

    The main theorem Theorem \ref{main} will follow if we show that
    the homomorphism in Proposition \ref{wdsurj} is an isomorphism.
    It therefore suffices to show the injectivity of
    $\bar{\Psi_{\lambda}}$.

    For this we prove the following lemma.

\blem\label{generators} Let
    $\mathcal{R}:=:=R(T^l)[x_1,\ldots, x_n]$. The ring
    $\mathcal{R}/\mathcal{I}_{\lambda}$ is generated as an
    $R(T^l)$-module by $n \choose \lambda$ elements.  \elem

    \begin{proof}
      Let
      $\mathcal{R}':=\mathbb{Z}[u_1,\ldots, u_l, x_1,\ldots, x_n]$.
      Let $\mathcal{I}_{\lambda}'$ be the ideal in the polynomial ring
      $\mathcal{R}'$ with its natural grading, generated by the
      following homogeneous elements of degree $d$
\[\sum_{0\leq k\leq
    d}(-1)^{d-k}e_{k}(x_{i_1}, x_{i_2},\ldots, x_{i_s})
  h_{d-k}({u_{\phi_{\lambda}(1)}},\ldots,
  {u_{\phi_{\lambda}(s+1-d)}})\] for $1\leq s\leq n$,
$1\leq i_1<\cdots<i_s\leq n$ and $d\geq s+1-q$ where
$q:=p_{\lambda^\vee}(s)$.

It siffices to find polynomials $\Phi_1(x),\ldots, \Phi_m(x)$ whose
classes generate $\mathcal{R}'/\mathcal{I}_{\lambda}'$ as a
$\mathbb{Z}[u_1,\ldots, u_l]$-module. Since $\mathcal{R}$
(resp. $R(T^l)$) is the localization of $\mathcal{R}'$
(resp. $\mathbb{Z}[u_1,\ldots, u_l]$) at $u_1\cdots u_l$, this will
imply that the classes of the polynomials
$\Phi_1(x),\ldots, \Phi_m(x)$ in $\mathcal{R}/\mathcal{I}_{\lambda}$
will generate $\mathcal{R}/\mathcal{I}_{\lambda}$ as an
$R(T^l)=\mathbb{Z}[u_1^{\pm 1},\ldots, u_l^{\pm 1}]$-module. The lemma
will then follow.

Now, by Theorem \ref{ahtheorem} we note that $\mathcal{R}'$ is
isomorphic to the ring
$\mathcal{S}/\mathcal{J}_{\lambda}\cong
H^*_{T^l}(\cf_{\lambda})$. Further, by \cite[Lemma 5.2]{ah} we know
that there exists polynomials $\Phi_1(x),\ldots, \Phi_m(x)$ which
generate $\mathcal{R}'/\mathcal{I}_{\lambda}'$ as a
$\mathbb{Z}[u_1,\ldots, u_l]$-module. Hence the lemma. \end{proof}

Now we prove Theorem \ref{main}.
\begin{proof} (Proof of the main theorem) We have shown in Theorem
  \ref{springercell} that $K_{T^l}(\cf_{\lambda})$ is a free
  $R(T^l)$-module of rank $n \choose \lambda$. Further, by Lemma
  \ref{generators} $\mathcal{R}/\mathcal{I}_{\lambda}$ is generated as
  an $R(T^l)$-module by $n\choose\lambda$ elements. Thus we will have
  a surjective $R(T^l)$-module homomorphism
     \[{\uppsi}:R(T^l)^{n\choose\lambda}\lra
       \mathcal{R}/\mathcal{I}_{\lambda}.\] Since by Proposition
     \ref{wdsurj}, $\bar{\Psi_{\lambda}}$ is surjective, we will get a
     surjective $R(T^l)$-module homomorphism
     $\bar{\Psi_{\lambda}}\circ \uppsi$ between two free of the same
     rank $n \choose \lambda$ over the intergral domain $R(T^l)$.
     This implies that $\bar{\Psi_{\lambda}}\circ \uppsi$ is an
     isomorphism and hence $\bar{\Psi_{\lambda}}$ is an isomorphism.
   \end{proof}

   \brem\label{familiarrelations}
   We have the identity
   \begin{align*}&\sum_{0\leq k\leq d}(-1)^{d-k}e_{k}(u_{\phi_{\lambda}(1)},\ldots,
     u_{\phi_{\lambda}(s)})h_{d-k}(u_{\phi_{\lambda}(1)},\ldots,
     u_{\phi_{\lambda}(s+1-d)})\\ &=e_d(u_{\phi_{\lambda}(s+2-d)},\ldots,
     u_{\phi_{\lambda}(s)})\\ &=0\end{align*} since the number of variables in the
   set $\{u_{\phi_{\lambda}(s+2-d)},\ldots, u_{\phi_{\lambda}(s)}\}$
   is less than $d$. Thus we can rewrite the relations defining
   $\mathcal{I}_{\lambda}$ in $\mathcal{R}$ as follows:

   \be\label{equivkthrelationscompact}\begin{array}{ll}\sum_{0\leq
       k\leq d}&(-1)^{d-k} \big[e_{k}(x_{i_1}, x_{i_2},\ldots,
     x_{i_s})-e_{k}(u_{\phi_{\lambda}(1)},\ldots,
     u_{\phi_{\lambda}(s)})\big]\cdot\\ 
     & h_{d-k}(u_{\phi_{\lambda}(1)},\ldots,
     u_{\phi_{\lambda}(s+1-d)})\end{array}\ee for $1\leq s\leq n$,
   $1\leq i_1<\cdots<i_s\leq n$ and $d\geq s+1-q$ where
   $q=p_{\lambda^{\vee}}(s)$. In particular, when
   $\lambda=(1,\ldots,1)$ and $\mathcal{F}_{\lambda}=\mathcal{F}$ then
   $\lambda^{\vee}=(n)$. Hence $p_{\lambda^{\vee}}(s)=0$ for
   $1\leq s<n$ and $p_{\lambda^{\vee}}(n)=n$. Moreover,
   $\phi_{\lambda}$ is the identity map. In this case the relation
   \eqref{equivkthrelationscompact} reduces to
   \be\label{flagpres1}\sum_{0\leq k\leq d}(-1)^{d-k}\big[e_{k}(x_{1},
   x_{2},\ldots, x_{n})-e_{k}(u_{1},\ldots, u_{n})\big] \cdot
   h_{d-k}(u_{1},\ldots, u_{n+1-d})\ee for $d\geq 1$.  The relations
   \eqref{flagpres1} can be readily seen to be equivalent to the
   relations \eqref{flagpres} given in Theorem \ref{eqflagpres}.
   \erem
\section{Relation with the ordinary $K$-ring}
Recall that we have the augmentation $\upepsilon :R(T^l)\lra \mathbb{Z}$
which sends the class $[V]$ of any $T^l$-representation $V$ to
$\mbox{dim}(V)$.  This gives $\mathbb{Z}$ a structure of
$R(T^l)$-module. Also we have a canonical $R(T^l)$-module structure on
$K_{T^l}(\mathcal{F}_{\lambda})$ through pull back via the structure
morphism or equivalently by the map which sends the class $[V]$ of a
$T^l$-representation to the trivial $T^l$-equivariant vector bundle on
$\mathcal{F}_{\lambda}\times V$. Similarly we have a
$\mathbb{Z}$-module structure on $K(\mathcal{F}_{\lambda})$ given by
pull back via the structure morphism or equivalently by the map
$\uptheta$ which sends any positive integer $n$ to the trivial vector
bundle over $\mathcal{F}_{\lambda}$ of rank $n$.  Let
${\mathrm f}:K_{T^l}(\mathcal{F}_{\lambda})\lra K(\mathcal{F}_{\lambda})$ denote
the forgetful homomorphism which sends the class of any
$T^l$-equivariant vector bundle to the class of the underlying vector
bundle forgetting the $T^l$-structure. (see
\cite{husemoller},\cite{segal}).
\bcor\label{equivariantformality} The map
\[{\mathrm F}:\mathbb{Z}\otimes_{R(T^l)} K_{T^l}(\mathcal{F}_{\lambda})\lra
  K(\mathcal{F}_{\lambda})\] induced by $\uptheta$ on the first factor
and by $\mathrm f$ on the second factor is an isomorphism of
$\mathbb{Z}$-modules. In other words the Springer variety
$\mathcal{F}_{\lambda}$ is weakly equivariantly formal for $K$-theory
\cite{hl}.\ecor {\bf Proof:} By Theorem \ref{springercell},
$K_{T^l}(\mathcal{F}_{\lambda})$ is a free $R(T^l)$-module of rank
$n\choose \lambda$. By \cite[Proposition 3.1]{su} we have that
$K(\mathcal{F}_{\lambda})$ is a free $\mathbb{Z}$-module of rank
$n\choose \lambda$ and also that it is generated by the class of the
line bundles $[L_i]$ for $1\leq i\leq n$. Since $L_i$ are
$T^l$-equivariant this implies that $[L_i]$ are the image of
$[L_i]_{T^l}\in K_{T^l}(\mathcal{F}_{\lambda})$ so that the map
$\mathrm F$ is
surjective. Thus $\mathrm F$ is a surjective homomorphism between two free
abelian groups of the same rank $n\choose \lambda$ and is therefore an
isomorphsm.  $\Box$

On the other hand by specializing to $u_i=1$ for all $1\leq i\leq l$,
that is by tensoring with $\mathbb{Z}$ over $R(T^l)$ we have
\[h_{d-k}({u_{\phi_{\lambda}(1)}},\ldots, {u_{\phi_{\lambda}(s+1-d)}})
  =\binom{s+1-d+d-k-1}{s+1-d-1}=\binom{s-k}{s-d}\]
(see \cite{mac}) for $d\geq s+1-q$.

Furthermore, since
\[\binom{q+d-k-1}{q-1}=\left(\prod_{i=1}^{q-s-1+d}\frac{s-k+i}{s-d+i}\right)\cdot
  \binom{s-k}{s-d}\] it follows that the equivariant $K$-theoretic
Tanisaki ideal $\mathcal{I}_{\lambda}$ reduces to the ordinary
$K$-theoretic Tanisaki ideal $I_{\lambda}$ defined in
\cite[p. 11]{su}.

By combining Theorem \ref{main} and Corollary
\ref{equivariantformality} we retrieve the presentation for the
ordinary $K$-ring of $\mathcal{F}_{\lambda}$ which was obtained
directly in \cite{su} by alternate methods.

\bth\label{ordk}[\cite[Theorem 4.2]{su}] The map
$\Psi^{\mbox{ord}}_{\lambda}:R=\mathbb{Z}[x_1,\ldots, x_n]\ra
K(\mathcal{F}_{\lambda})$ be the ring homomorphism defined by
$\Psi^{\mbox{ord}}_{\lambda}(x_j)=[L_j]$ for $1\leq j\leq n$. Then
$\Psi^{\mbox{ord}}_{\lambda}$ is surjective and
$ker(\Psi^{\mbox{ord}}_{\lambda})=I_{\lambda}$ where $I_{\lambda}$ is
the ideal in $R$ generated by the elements
\[ \sum_{0\leq k\leq d}(-1)^{d-k}e_k(x_{i_1},\ldots, x_{i_s})\cdot
  \binom{q+d-k-1}{q-1}\] where $1\leq
i_1<\cdots<i_s\leq n$, $1\leq s\leq n$ and $d\geq s+1-q$.
\eeth

\vspace{0.1cm} {\bf Acknowledgement:} The work was supported by Serb
Matrics research grant with project number MTR/2022/000484. The author
is very grateful to Prof. Parameswaran Sankaran for several valuable
discussions. I am also grateful to him for a careful reading of the earlier
versions of the manuscript and for sharing several valuable comments and
suggestions. The author wishes to thank Prof. Shrawan Kumar for very
valuable suggestions and encouragement during the preparation of the
manuscript. The author also wishes to thank Prof. Megumi Harada for
encouragement and helpful email exchanges. The author thanks the
referee for very valuable comments and suggestions.


\begin{thebibliography}{99}
  
\bibitem{ah} Abe, H. and Horiguchi, T., The torus equivariant
  cohomology rings of Springer varieties, Topology and its
  Applications {\bf 208}, 143-159 (2016)

\bibitem{at} Atiyah, M.F., K-Theory, W.A.Benjamin, New York, NY, (1967).

\bibitem{athir} Atiyah, M.F., and Hirzebruch, F., Vector bundles and
  homogeneous spaces, Proceedings of Symposia in Pure Mathematics,3,
  Amer. Math. Soc. (1961).

  \bibitem{atiyahtall} M.F. Atiyah and D. O. Tall, {\it Group
    Representations, A-rings and the J-homomorphism}, Topology {\bf 8}
  (1969), pp. 253-297.

\bibitem{bb} A. Bialynicki-Birula, {\it Some properties of the
    decomposition of algebraic varieties determined by the action of a
    torus}, Bull. Acad. Polon. Sci. Math. Astronom. Phys. {\bf 24}
  (1976), 667-674.
  
\bibitem{borel} Borel, A., Sur La cohomologie des espaces fibres principaux et des espaces homogenes de groupes de Lie compacts. Ann. Math. {\bf 57} (1953), 115--207.

  
\bibitem{dp} De Concini, C. and Procesi, C.,  Symmetric
  Functions, Conjugacy Classes and the Flag Variety,
  Invent. Math. {\bf 64}, 203-219 (1981).


\bibitem{hl} Megumi Harada and Gregory D. Landweber, Surjectivity for
  Hamiltonian $G$-Spaces in K-Theory {\it Transactions of the
    American Mathematical Society} {\bf 359}, No. 12, 6001-6025
  (2007).

  
\bibitem{hs} Hotta, R. and Springer T. A., A Specialization theorem for
  certain Weyl group representations and an application to Green
  polynomials of unitary groups, Invent. Math. {\bf 41}, 113-127 (1977).

\bibitem{husemoller} Husemoller, D., {\it Fibre bundles} 2nd ed. GTM{\bf 20}  Springer-Verlag, New York, 1975. 


\bibitem{k} Karoubi, M., {\it K-theory}, Grund. Math. Wiss.226,
  Springer-Verlag, Berlin, (1978).

\bibitem{kk} B. Kostant and S. Kumar, {\it $T$-equivariant $K$-theory
    of generalized flag varieties} J. Differential Geom. {\bf 32}
  (1990),549-603.


\bibitem{mac} Macdonald, I. G., {\it Symmetric functions and Hall
    polynomials}, Oxford Univ. Press 1979.

\bibitem{mcleod} J. McLeod, {\it The Kunneth formula in equivariant
    $K$-theory}, in {\it Algebraic Topology, Waterloo,} 1978
  (Proc. Conf. Univ. Waterloo, Waterlo, Ont., 1978), Lecture Notes in
  Mathematics, Vol. 741, Springer-Verlag, Berlin, 1979, pp. 316-333. 

\bibitem{su} Sankaran, P., and Uma, V., K-theory of Springer
  varieties, math arXiv:2201.03058, to appear in Tohoku Math. Journ.

  
\bibitem{segal} Segal, G., Equivariant K-theory, {\it  Publications
  mathématiques de l’I.H.É.S.}, tome 34 (1968), p. 129-151

\bibitem{spal} Spaltenstein, N., On the fixed point set of a unipotent
  transformation on the flag manifold,
  Nederl. Akad. Wetensch. Proc. Ser. A {\bf 79} (1976) 452-456.

\bibitem{spr1} Springer, T. A., Trignometric sums, Green functions of
  finite groups and representations of Weyl groups, Invent. Math. {\bf
    36} (1976), 173-207

\bibitem{spr2} Springer, T. A., A construction of representations of Weyl
  groups, Invent. Math. {\bf 44} (1978) 279-293.


\bibitem{t} Tanisaki, T., Defining Ideals of the Closures of the
  Conjugacy Classes and Representations of the Weyl groups, Tohoku
  Math. Journ.{ \bf 34} (1982) 575-585.

\bibitem{tym} Tymoczko, J., The Geometry and Combinatorics of
  Springer Fibers, arXiv:1606.02760v1[math.AG].  











\end{thebibliography}
\end{document}